\documentclass{article}

\usepackage{epsfig}
\usepackage{amsmath, amssymb, latexsym}
\usepackage{graphs}

\newtheorem{theorem}{Theorem}[section]
\newtheorem{proposition}[theorem]{Proposition}
\newtheorem{lemma}[theorem]{Lemma}
\newtheorem{corollary}[theorem]{Corollary}

\newcommand{\N}{{\mathbb N}}
\newcommand{\Z}{{\mathbb Z}}
\newcommand{\R}{{\mathbb R}}
\newcommand{\E}{{\mathcal E}}
\def\es{{\emptyset}}

\def\Spec{{\rm Spec}}
\def\ove{\ov{e}}
\def\ov{\overline}
\newcommand{\lam}{\lambda}
\def\Lam{\Lambda}
\newcommand{\Ek}{{\mathcal E}}
\newcommand{\om}{\omega}
\def\ovom{\ov{\omega}}
\def\nutil{\wtil{\nu}}
\def\wtil{\widetilde}
\def\Rk{{\mathcal R}}
\def\Dk{{\mathcal D}}
\def\ovtau{\ov{\tau}}
\newcommand{\gam}{\gamma}
\newcommand{\eps}{{\varepsilon}}
\def\th{\theta}
\def\bx{{\bf x}}
\def\dist{{\rm dist}}
\newcommand{\Q}{{\mathbb Q}}
\def\qed{\hfill$\square$}

\def\sig{\sigma}
\def\Lk{{\mathcal L}}
\def\Ak{{\mathcal A}}
\def\Ik{{\mathcal I}}
\def\Jk{{\mathcal J}}
\def\Ok{{\mathcal O}}
\newcommand{\countnow}{\refstepcounter{theorem}{\bf\arabic{section}.\arabic{theorem}.} }
\newcommand{\mydefinition}{\medskip{\bf Definition }\countnow}
\newcommand{\myremark}{\medskip{\bf Remark }\countnow}
\newcommand{\myexample}{\medskip{\bf Example }\countnow}

\begin{document}


\title{Invariant Measures on Stationary Bratteli Diagrams}

\author{{\bf S.~Bezuglyi}\\
Institute for Low Temperature Physics, Kharkov, Ukraine \\
bezuglyi@ilt.kharkov.ua\\
 {\bf J.~Kwiatkowski}\footnote{The research of J.K was supported by  grant MNiSzW  nr NN201384834.}\\
 College of Economics and Computer Sciences, Olsztyn, Poland \\
 jkwiat@mat.uni.torun.pl
 \\
 {\bf K. Medynets}\footnote{K.M. was supported by the Akhiezer fund and INTAS YSF 05-109-53-15.}
 \\
 Institute for Low Temperature Physics,  Kharkov, Ukraine
 \\
 medynets@ilt.kharkov.ua
 \\
 {\bf B. Solomyak} \footnote{ B.S. was supported in part by NSF grants DMS-0355187 and DMS-0654408.}\\
University of Washington, Seattle, USA
\\
solomyak@math.washington.edu
 }
 
 \date{}
   
  \maketitle


\begin{abstract}
We study dynamical systems acting on the path space of
a stationary (non-simple) Bratteli diagram. For such systems we
explicitly describe all ergodic probability measures  invariant with
respect to the tail equivalence relation (or the Vershik map). These
measures are completely described by the incidence matrix of the
diagram. Since such diagrams correspond to substitution dynamical
systems,
 this description gives an algorithm for finding invariant probability  measures
 for aperiodic non-minimal substitution systems. Several corollaries of these results are obtained.
 In particular, we show that the invariant measures are not mixing and give a criterion for
 a complex number to be an eigenvalue for the Vershik map.
\end{abstract}

\section{Introduction}
Every homeomorphism $T$ of a compact metric space has a nontrivial
set of $T$-invariant Borel probability measures. This set forms a
simplex in the set of all probability invariant measures whose
extreme points are ergodic $T$-invariant measures. There is an extensive list
of research papers devoted to the study of  relations between
properties of transformations and those of the corresponding simplex
of invariant measures. We mention only some relatively recent papers
by Akin \cite{akin:1999,akin:2004}, Downarowicz
\cite{Downarowicz:1991,Downarowicz:2006}, Glasner and Weiss,
\cite{glasner_weiss:1995,glasner_weiss:1997}, Gjerde and Johansen \cite{GJ}, a few older ones \cite{bauer_sigmund:1975,sigmund:1978}, and the
well-known books on ergodic theory \cite{walters:book},
\cite{petersen:book}, \cite{csf:book}.
Any aperiodic transformation in measurable,
Borel, and Cantor dynamics can be realized as a Vershik map acting
on the path space of a Bratteli diagram \cite{vershik:1981,vershik:1982},
\cite{herman_putnam_skau:1992},
\cite{bezuglyi_dooley_kwiatkowski:2006-1}, \cite{medynets:2006}.
Such a representation of aperiodic transformations is very
convenient from various  viewpoints, in particular, for finding
invariant measures and their values on clopen sets. We should note
here that the converse statement is not, in general, true in the
framework of Cantor dynamics: there are Bratteli diagrams which do
not admit continuous Vershik maps \cite{medynets:2006}. The
suggested  approach naturally leads us to study  probability
measures on the path spaces of Bratteli diagrams which are invariant
with respect to the tail (cofinal) equivalence relation.
Such measures also arise as {\it states of the dimension group}
associated with the Bratteli diagram, see \cite{Effros}. They were considered by
Kerov and Vershik \cite{kerov_vershik}, who called them
{\it central measures} since they appeared as central states on certain
$C^*$-algebras. There are some classes of Bratteli diagrams for which the invariant measures are known, but the focus has been either on uniquely ergodic systems,
e.g.\ simple stationary diagrams \cite{durand_host_skau}, linearly recurrent systems \cite{CDHM}, or very specific cases, such as
the Pascal diagram \cite{petersen_schmidt} or Euler diagram \cite{BKPS}.
Non-simple stationary diagrams have not been studied systematically.

The main goal of the present paper is to give an explicit
description of  probability measures on the path space of a
stationary Bratteli diagram which are invariant with respect to the
tail equivalence relation, assuming that this equivalence relation
is aperiodic. We describe our main results briefly (precise definitions and statements are given later). A stationary Bratteli diagram is  determined by its incidence matrix $F$.
It is well-known that for {\it simple} stationary Bratteli diagrams, i.e.\ when the incidence matrix is
primitive, the invariant probability measure is unique and determined by the Perron-Frobenius (PF)
eigenvector of $A=F^T$. (This is proved in Effros \cite[Theorem 6.1]{Effros} using the language of states and dimension groups.
Fisher \cite{Fisher1} points out that this result is implicitly contained in
\cite[Lemma 2.4]{Bowen_Marcus}.) In the general case, we prove that finite invariant measures are in 1-to-1
correspondence with the {\it core} of $A$, defined by $core(A) = \bigcap_{n=0}^\infty A^n(\R_+^N)$ where
$A$ has size $N\times N$. Perron-Frobenius theory for non-negative matrices (see \cite{schneider})
says that $core(A)$ is a simplicial cone, and when the irreducible components of $A$ are primitive
(which can always be achieved by ``telescoping''), its extremal rays are generated by non-negative
eigenvectors of $A$. Every such an eigenvector is the PF eigenvector for one of the irreducible components
of $A$, but only the  {\em distinguished} components yield a non-negative eigenvector.
A component $\alpha$ is distinguished if its PF eigenvalue is strictly greater than PF eigenvalues of all components which have access to
 $\alpha$, see Section 3 for definitions.
Thus, ergodic invariant probability measures are in 1-to-1 correspondence with distinguished components
(Theorem~\ref{Theorem_MeasuresErgodicMeasures}). Interestingly, non-distinguished components also play a role; in fact, they are in 1-to-1
correspondence, up to a constant multiple,
with ergodic $\sigma$-finite (infinite) measures that are positive and finite on some open set (Theorem~\ref{th-infmeas}).
We should note that some of our results are implicitly contained in \cite{Handelman}, see Remark 3.8.3.

Substitution dynamical systems have been studied extensively; however, in the vast majority of papers,
primitivity, hence minimality, is assumed. Every primitive substitution system is conjugate to the
Vershik map on a simple stationary Bratteli diagram \cite{forrest,durand_host_skau}, and this has
recently been extended to a large class of aperiodic non-minimal substitutions in
\cite{bezuglyi_kwiatkowski_medynets:substitutions}.
Thus, our results yield an explicit description of invariant
measures (both finite and $\sigma$-finite) for such systems.
In contrast to the case of minimal substitution systems (see
\cite{queffelec:book}), aperiodic substitution systems are not, in
general, uniquely ergodic. For instance, consider the
following two substitution systems $(X_\sigma,T_\sigma)$ and
$(X_\tau,T_\tau)$ defined on the alphabet $A =\{a,b,c\}$ by
substitutions $\sigma$ and $\tau$ where  $\sigma(a)=\tau(a)=abb$,
$\sigma(b)=\tau(b)=ab$, and $\sigma(c)=accb$, $\tau(c)=acccb$. Each
of these systems has a unique minimal component $C$. However, it
follows from our results that
$(X_\sigma,T_\sigma)$ has a unique invariant probability measure supported on
the minimal component, whereas $(X_\tau,T_\tau)$ has two
ergodic invariant probability measures: one of them is supported on $C$ and the
other measure is supported on the complement of $C$. Thus, these
systems cannot be conjugate and they cannot even be orbit
equivalent.

Recently Yuasa
\cite{yuasa} obtained a somewhat similar result for ``almost-minimal'' substitutions, which is
complementary to ours, since those substitution systems have a fixed point (for the shift transformation). Earlier, a special
case of such a substitution, namely, $0\to 000, \ 1\to 101$,  was studied by Fisher \cite{fisher}.

The set of invariant measures is of crucial importance for the
classification of Cantor minimal systems up to orbit equivalence
\cite{giordano_putnam_skau:1995,glasner_weiss:1995-1}.
Recall the following results proved by Giordano, Putnam, and Skau in
\cite{giordano_putnam_skau:1995}: (1) two Cantor minimal systems
$(X,T)$ and $(Y,S)$   are orbit equivalent if and only if there
exists a homeomorphism $F : X\to Y$ carrying the $T$-invariant
probability measures onto the $S$-invariant probability measures;
(2) two uniquely ergodic Cantor minimal systems $(X, T)$ and $(Y,
S)$ are orbit equivalent if and only if the clopen values sets for $\mu$ and $\nu$ are the same, i.e.\ $\{\mu(E) : E \mbox{ clopen
in } X\} = \{\nu(F) : F \mbox{ clopen in } Y\}$ where $\mu$ and
$\nu$ are unique probability invariant measures for $T$ and $S$
respectively. The notion of orbit equivalence for {\it aperiodic}
Cantor systems has not been studied yet.  Based on our study of stationary
Bratteli diagrams,  we show that the second statement does not hold
any more for non-minimal uniquely ergodic homeomorphisms.
We intend to apply our results to the study of orbit equivalence of aperiodic homeomorphisms of a Cantor set in another paper.

We also study some properties of measure-preserving systems on stationary diagrams, corresponding to the ergodic probability measures.
In particular, we show that they are not mixing and give a
criterion for a complex number to be an eigenvalue.
These results have common features with some in the literature, see e.g. \cite{DK,Liv88,CDHM,BDM}, but they do not follow from them, since minimality, and hence unique ergodicity has been a common assumption until now.

The article is organized as follows. Section 2 contains some
definitions and facts concerning Bratteli diagrams which are used in
the subsequent sections. We also discuss the construction of
invariant measures on Bratteli diagrams of general form. Section 3
is focused on the proof of the main result which gives an explicit
description of invariant probability measures. In Section 4 we obtain further properties of these measures and describe $\sigma$-finite
invariant measures.  Section 5 contains several applications of our
results and examples. In particular, we show
that ergodic invariant probability measures  of a substitution
aperiodic system can be determined from eigenvectors of the
incidence matrix of the substitution. In Section 6, we study ergodic-theoretic properties of our systems.

%
%

\section{Measures on Bratteli diagrams}\label{section2} In this section, we study  Borel  measures on the path space of a Bratteli diagram which are invariant with respect to the tail equivalence relation.  Since the notion of Bratteli diagrams has
been discussed in many well-known  papers on Cantor dynamics (e.g.
\cite{herman_putnam_skau:1992} and
\cite{giordano_putnam_skau:1995}), we present here the main
definitions and notation only.  We also refer the reader to the
works   \cite{medynets:2006} and
\cite{bezuglyi_kwiatkowski_medynets:substitutions}, where
Bratteli-Vershik models of Cantor aperiodic systems and aperiodic
substitution systems were considered.

\mydefinition \label{Definition_Bratteli_Diagram} A {\it Bratteli
diagram} is an infinite graph $B=(V,E)$ such that the vertex set
$V=\bigcup_{i\geq 0}V_i$ and the edge set $E=\bigcup_{i\geq 1}E_i$
are partitioned into disjoint subsets $V_i$ and $E_i$ such that

(i) $V_0=\{v_0\}$ is a single point;

(ii) $V_i$ and $E_i$ are finite sets;

(iii) there exist a range map $r$ and a source map $s$ from $E$ to
$V$ such that $r(E_i)= V_i$, $s(E_i)= V_{i-1}$, and
$s^{-1}(v)\neq\emptyset$, $r^{-1}(v')\neq\emptyset$ for all $v\in V$
and $v'\in V\setminus V_0$.
 \medbreak

The pair  $(V_i,E_i)$ is called the $i$-th level of the diagram $B$.
We write $e(v,v')$ to denote an edge $e$ such that $s(e)=v$ and
$r(e)=v'$.

A finite or infinite sequence of edges $(e_i : e_i\in E_i)$ such
that $r(e_{i})=s(e_{i+1})$ is called a {\it finite} or {\it infinite
path}, respectively.  For a Bratteli diagram $B$, we denote by $X_B$
the set of infinite paths starting at the  vertex $v_0$. We endow
$X_B$ with the topology generated by  cylinder sets
$U(x_1,\ldots,x_n):=\{x\in X_B : x_i=e_i,\;i=1,\ldots,n\}$, where
$(e_1,\ldots,e_n)$ is a finite path in $B$. Then $X_B$ is a
0-dimensional compact metric space with respect to this topology. We
will consider such diagrams $B$ for which the path space $X_B$ has
no isolated points.

Each Bratteli diagram can be given a diagrammatic representation
(see, for instance, Fig. 1).

\unitlength = 0.5cm
\begin{center}
\begin{graph}(10,13)
\graphnodesize{0.4}
%
\roundnode{V0}(5,12)
\roundnode{V11}(1,9) \roundnode{V12}(5,9) \roundnode{V13}(9,9)
\edge{V11}{V0} \edge{V12}{V0} \edge{V13}{V0}
\roundnode{V21}(1,4.5) \roundnode{V22}(5,4.5) \roundnode{V23}(9,4.5)
\bow{V21}{V11}{0.06} \bow{V21}{V11}{-0.06} \bow{V21}{V12}{0.06}
\bow{V21}{V12}{-0.06} \bow{V22}{V12}{0.06} \bow{V22}{V12}{-0.06}
\edge{V22}{V11} \bow{V23}{V13}{0.06} \bow{V23}{V13}{-0.06}
\edge{V23}{V11} \edge{V23}{V12}
\roundnode{V31}(1,0.5) \roundnode{V32}(5,0.5) \roundnode{V33}(9,0.5)
\bow{V31}{V21}{0.06} \bow{V31}{V21}{-0.06} \bow{V31}{V22}{0.06}
\bow{V31}{V22}{-0.06} \bow{V32}{V22}{0.06} \bow{V32}{V22}{-0.06}
\edge{V32}{V21} \bow{V33}{V23}{0.06} \bow{V33}{V23}{-0.06}
\edge{V33}{V21} \edge{V33}{V22}
\end{graph}
\vskip0.2cm $.\ .\ .\ .\ .\ .\ .\ .\ .\ .\ .\ .\ .\ .\ .\ .$
\vskip0.3cm Fig. 1
\end{center}

Given a Bratteli diagram $B=(V,E)$, fix a level $n\geq 1$. Define the
$|V_{n+1}|\times |V_n|$ matrix $F_{n}=(f^{(n)}_{vw})$ whose  entries
$f^{(n)}_{vw}$  are equal to  the number of edges between the
vertices $v\in V_{n+1}$ and $w\in V_{n}$, i.e.,
$$
 f^{(n)}_{vw} = |\{e\in E_{n+1} : r(e) = v, s(e) = w\}|.
$$
(Here and thereafter $|A|$ denotes the cardinality of the set $A$.)
For instance, we have that for the above diagram
$$F_{1}= F_2 =\left(\begin{array}{ccc} 2 & 2 & 0\\
1 & 2 & 0\\
1 & 1 & 2 \end{array}\right).
$$

A Bratteli diagram $B=(V,E)$ is called {\it stationary} if
$F_{n}=F_1$ for every $n\geq 2$.

Observe that every vertex $v\in V $ is connected to $v_0$ by a
finite path, and the set $E(v_0,v)$ of all such paths  is finite.
Set  $h_v^{(n)}=|E(v_0,v)|$ where $v\in V_{n}$ and $h^{(n)}=(h_w^{(n)})_{w\in V_n}$.
Then we get that for all $n\ge 1$,
\begin{equation}\label{heights}
h^{(n+1)}=\sum_{w\in V_{n}}f_{vw}^{(n)}h^{(n)}_w=F_{n}h^{(n)}.
\end{equation}

For $w\in V_n$, the set $E(v_0,w)$ defines the clopen subset
$$
X^{(n)}_w=\{x=(x_i)\in X_B : r(x_n)=w \}.
$$
Moreover, the sets  $\{X^{(n)}_w : w\in V_n\}$ form a clopen
partition of $X_B,\ n\geq 1$. Analogously, each finite path
$\overline e=(e_1,\ldots,e_n)\in E(v_0,w)$ determines the clopen set
$$
X^{(n)}_w(\overline e)=\{x=(x_i)\in X_B : x_i=e_i,\;i=1,\ldots,n \}.
$$
These sets form a clopen partition of $X_w^{(n)}$. We will also
use the notation $[\overline{e}]$ for the clopen set
$X^{(n)}_w(\overline e)$ if it does not lead to a confusion.
\medbreak

By definition, a Bratteli diagram $B =(V,E)$ is called {\it ordered}
if every set $r^{-1}(v)$, $v\in \bigcup_{n\ge 1} V_n$,   is linearly
ordered, see \cite{herman_putnam_skau:1992}. Denote by $\mathcal O =
\mathcal O(B) $ the set of all possible orderings on $B$.  An
ordered Bratteli diagram will be denoted by $B(\omega) = (V,E,
\omega)$ where $\omega \in \mathcal O$. Given $B(\omega)$,  any two
paths from $E(v_0,v)$ are comparable with respect to the
lexicographical order. We call a
finite or infinite path $e=(e_i)$ {\it maximal (minimal)} if every
$e_i$ is maximal (minimal) amongst the edges from $r^{-1}(r(e_i))$.
Notice that for  $v\in V_i,\ i\ge 0$, the minimal and maximal
(finite) paths in $E(v_0,v)$ are unique. Denote by
$X_{\max}(\omega)$ and $X_{\min}(\omega)$ the sets of all maximal
and minimal infinite paths from $X_B$, respectively. It is not hard
to see that $X_{\max}(\omega)$ and $X_{\min}(\omega)$ are non-empty
closed subsets.

A Bratteli diagram $B=(V,E, \omega)$ is called {\it stationary
ordered} \cite{durand_host_skau} if it is stationary and the partial linear order on $E_n$,
defined by $\omega$, does not depend on $n$.

Let $B = (V,E, \omega)$ be a stationary ordered Bratteli diagram, or more generally, suppose $N=\sup_n |V_n| < \infty$.
Then it is easy to see that the sets $X_{\max}(\omega)$ and $X_{\min}(\omega)$ of maximal
and minimal paths are finite
\cite{bezuglyi_kwiatkowski_medynets:substitutions}. Indeed, observe that two maximal paths which go through the same vertex at level $n$ must have the same beginning $e_1,\ldots,e_n$. Given $N+1$ maximal paths, we can find two of them which go through the same vertex at infinitely many levels, hence they must coincide. 

\mydefinition Let $B = (V,E,\omega)$ be an  ordered Bratteli
diagram. We say that $\varphi = \varphi_\omega : X_B\rightarrow X_B$
is a {\it Vershik map} if it satisfies the following conditions:

(i) $\varphi$ is a homeomorphism of the Cantor set $X_B$;

(ii) $\varphi(X_{\max}(\omega))=X_{\min}(\omega)$;

(iii) if an infinite path $x=(x_1,x_2,\ldots)$ is not in
$X_{\max}(\omega)$, then
$\varphi(x_1,x_2,\ldots)=(x_1^0,\ldots,x_{k-1}^0,\overline
{x_k},x_{k+1},x_{k+2},\ldots)$, where $k=\min\{n\geq 1 : x_n\mbox{
is not maximal}\}$, $\overline{x_k}$ is the successor of $x_k$ in
$r^{-1}(r(x_k))$, and $(x_1^0,\ldots,x_{k-1}^0)$ is the minimal path
in $E(v_0,s(\overline{x_k}))$.

If $f$ is a Borel automorphism of $X_B$ which satisfies conditions
(ii) and (iii), then  $f$ is called a Borel-Vershik automorphism.
\medbreak
\myremark 1. Vershik maps were introduced in \cite{vershik:1981,vershik:1982} in the measure-theoretic category, where they are called
{\em adic transformations}.

2. A Vershik map acts as the ``immediate successor transformation'' in the (reverse) lexicographic ordering induced by $\om$ on
$X_B \setminus X_{\max}$, and it is easily seen to be continuous on this set. In order to get a homeomorphism, one needs to map
$X_{\max}$ onto $X_{\min}$ bijectively, and of course, to have continuity of this extension and its inverse.
It is shown in \cite{medynets:2006} that there are stationary
Bratteli diagrams which do not admit a Vershik map.

\mydefinition Let $B=(V,E)$ be a Bratteli diagram.
Two infinite paths $x=(x_i)$ and $y=(y_i)$ from $X_B$ are said to be {\it tail
equivalent} if    there exists $i_0$ such that $x_i = y_i$ for all
$i\geq i_0$. Denote by $\mathcal R$  the tail equivalence relation
on $X_B$. 

\medbreak

We mention here the work \cite{giordano_putnam_skau:2004}
where  various properties of the tail equivalence relations  are
discussed in the context of Cantor dynamics.  \medbreak

\mydefinition A Borel equivalence relation is called {\it aperiodic} if all the equivalence classes are infinite.

\medbreak

{\it Throughout the paper, we consider Bratteli diagrams $B$ for
which $\mathcal R$ is an aperiodic Borel equivalence relation on
$X_B$.} In other words, every $\Rk$-equivalence class is countably infinite (it is obviously at most countable).

\myremark
Observe that the Vershik map (if it exists) is uniquely determined by the order $\om \in \Ok$ if the set $X_{\max}(\om)$ has empty
interior. 
One can show that $int(X_{\max}(\omega)) \neq \emptyset$ (or
$int(X_{\min}(\omega)) \neq \emptyset$) if and only if there exist
$n_0 \in \N$ and $x = (x_i)\in  X_B$ such that the cylinder set $U(x_1,\ldots,x_n)=
\{y = (y_i) \in X_B : y_1 = x_1,\ldots y_n =
x_n\}$ has no distinct cofinal paths for all $n> n_0$. It follows that $int(X_{\max}(\omega))=\es$ if and only if the equivalence relation $\Rk$ is aperiodic.

\medskip

For a Bratteli diagram $B$, denote by $M(\Rk)$ the set of finite positive Borel $\Rk$-invariant measures, and by
$M_1(\mathcal R)\subset M(\Rk)$ the set of invariant probability measures.
Similarly, $M_\infty(\mathcal R)$ denotes the set of non-atomic $\sigma$-finite
infinite $\mathcal R$-invariant measures. (We will
use below the term ``infinite measure'' for a $\sigma$-finite
infinite non-atomic measure.) Recall that a measure $\mu$ is called
{\it $\mathcal R$-invariant} if it is invariant under the Borel
action of any countable group $G$ on $X_B$ whose orbits generate the
equivalence relation $\mathcal R$. For a Bratteli diagram,  such a
group $G$ can be chosen  locally finite; it is sometimes called the group of ``finite coordinate changes.''

A Borel measure on $X_B$ is completely determined by its values on cylinder sets, since they generate the Borel $\sigma$-algebra.
Thus, we have that $\mu$ is $\Rk$-invariant if and only if for any $n$ and any $w\in V_n$,
\begin{equation} \label{eq-invar}
\overline{e}, \overline{e}' \in E(v_0, w) \ \Longrightarrow\
\mu(X_w^{(n)}(\overline e)) =\mu(X_w^{(n)}(\overline e')).
\end{equation}

\begin{lemma}\label{invariant measures} Let $B=(V,E, \omega)$ be an ordered Bratteli diagram which admits an aperiodic Vershik map $\varphi_\omega$, and  suppose that
the tail equivalence relation $\mathcal R$ is aperiodic. Then the set $M_1(\mathcal R)$ coincides with the set $M_1(\varphi_\omega)$ of $\varphi_\omega$-invariant probability measures. Furthermore, $M_\infty(\mathcal R) = M_\infty(\varphi_\omega)$ for a stationary Bratteli diagram $B$.
\end{lemma}

{\it Proof}. If
$x\in X\setminus Orb_\varphi(X_{\max}(\omega)\cup X_{\min}(\omega))$, then the $\varphi_\omega$-orbit
of $x$ is equal to the equivalence class $\mathcal R(x)$. If $\mu$
is a  $\varphi_\omega$-invariant finite measure, then
$\mu(X_{\min}(\omega))= \mu(X_{\max}(\omega)) =0$ because the sets 
$X_{\min}(\omega)$ and $X_{\max}(\omega)$ are wandering with respect to $\varphi_\omega$ (we are using here the aperiodicity assumption).
The proof of the relation $M_\infty(\mathcal R) =
M_\infty(\varphi_\omega)$ for stationary diagrams follows from the fact that
$X_{\min}(\omega)$ and $X_{\max}(\omega)$ are finite sets.
\hfill$\square$
\\

It follows from this lemma that for an ordered Bratteli diagram $B =
(V,E,\omega)$ and an $\mathcal R$-invariant measure $\mu$ we can
study properties of the measure-theoretical   dynamical system $(X_B,
\mu, \varphi_\omega)$ independently of whether the
Vershik map $\varphi_\omega$ exists everywhere on $X_B$.
\medskip

Let $B = (V,E,\omega)$ be an ordered Bratteli diagram. It is
clear that for the sets of infinite invariant measures we have the
relation $M_\infty(\mathcal R)\supseteq M_\infty(\varphi_\omega)$.
We do not know whether these sets are always equal. It would be so if
we could show that $\mu(X_{\min}(\omega))= 0$ for any infinite
$\varphi_\omega$-invariant non-atomic measure $\mu$. \medbreak

Let us consider next the case of finite $\mathcal R$-invariant
measures for a Bratteli diagram $B=(V,E)$. Take a Borel
measure $\mu \in M(\mathcal R)$. Recall that such a measure is
uniquely determined by its values on clopen sets of $X_B$. This
means that if we know $\mu(X^{(n)}_w)$ for all $w\in V_n$ and $n
\geq 1$, then $\mu$ is completely defined. In view of (\ref{eq-invar}),
$$
\mu(X^{(n)}_w(\overline e)) = \frac{1}{h_w^{(n)}}\mu(X^{(n)}_w),\
\overline e \in E(v_0,w).
$$
Set $p^{(n)} = (p_w^{(n)})_{w\in V_n}$ where $p_w^{(n)} =
\mu(X^{(n)}_w(\ov{e})),\ n\geq 1$, for some $\ov{e}\in E(v_0,w)$.
For $\overline e\in E(v_0,w)$ and $e\in E(w,v)$, $v\in V_{n+1}$,
denote by $(\overline ee)$ the finite path that coincides with
$\overline e$ on first $n$ segments and whose $(n+1)$-st edge is
$e$. Thus, we get a disjoint union
$$X_w^{(n)}(\overline e)=\bigcup_{v\in V_{n+1}}\bigcup_{e\in E(w,v)}X_v^{(n+1)}(\overline ee).$$
It follows that
\begin{eqnarray*}
\mu(X_w^{(n)}(\overline e))&=&\sum_{v\in V_{n+1}}\sum_{e\in
E(w,v)}\mu(X_v^{(n+1)}(\overline ee))\\& =& \sum_{v\in V_{n+1}}
f^{(n)}_{vw}\mu(X_v^{(n+1)}(\overline ee))\\&=&\sum_{v\in V_{n+1}}
f^{(n)}_{vw}p_v^{(n+1)}
\end{eqnarray*}
where $f^{(n)}_{vw}$ are the entries of $F_n$. Thus,
\begin{equation} \label{from}
p^{(n)} = F_n^T p^{(n+1)},\ n\ge 1.
\end{equation}

We recall standard definitions pertaining to cones.

\mydefinition \label{def-cone}
A subset $C\subset \mathbb R^N$ is called a {\it convex cone} if $\alpha
x+\beta y\in C$ for all $x,y\in C$ and $\alpha,\beta \geq 0$. A
subcone $Q$ of $C$ is called a {\it face} of $C$ if $x\in Q$, $y\in
C$, and $x-y\in C$ imply $y\in Q$. For $x\in C$, denote by $\Phi(x)$
the minimal face (the  intersection of all faces) that contains $x$.
A vector $x\in C$ is called {\it an extreme vector} if $\Phi(x)$ is
the ray generated by $x$, i.e. $\Phi(x)=\{\alpha x : \alpha\geq
0\}$. In this case, $\Phi(x)$ is also called an {\it extreme ray}.
A cone $C$ is called {\it polyhedral} ({\it finitely
generated}) if it has finitely many extreme rays. The cone is
{\it simplicial} if it has exactly $m$ extreme rays, where
$m=\mbox{dim(span }C)$. \medbreak

Now we go back to the context and notation of a Bratteli diagram.
For $x=(x_1,\ldots,x_N)\in \R^N$, we will write $x\ge 0$ if $x_i\ge 0$ for all $i$, and consider the positive cone
$\R^N_+=\{x\in\R^N : x\geq 0\}$. Let
$$
C_k^{(n)}:= F_k^T\cdots F_n^T\left(\R_+^{|V_{n+1}|}\right),\ \ 1 \le k \le n.
$$
Clearly, $\R_+^{|V_k|}\supset C_k^{(n)} \supset C_k^{(n+1)}$ for all $n\ge 1$. Let
$$
C_k^{\infty} = \bigcap_{n\ge k} C_k^{(n)}, \ k\ge 1.
$$
Observe that $C_k^\infty$ is a closed non-empty convex subcone of $\R_+^{|V_k|}$.
It also follows from these definitions that
$$
    F_k^T C_{k+1}^\infty = C_k^{\infty}.
$$
In general, the cones $C_k^{\infty}$ need not be simplicial, since there exist Bratteli diagrams with infinitely many ergodic invariant
measures. However, we  show in the next section that for stationary Bratteli
diagrams they are always simplicial.

The following result is formulated for finite $\mathcal R$-invariant
measures. The case of infinite measures is discussed in Remark
\ref{InfiniteMeasure}.

\begin{theorem}\label{Theorem_measures_general_case}  Let $B = (V,E)$ be a Bratteli diagram such that the tail equivalence relation
 $\mathcal R$ on $X_B$ is aperiodic.
If  $\mu\in M(\Rk)$,
then the vectors $p^{(n)}=(\mu(X_w^{(n)}(\ov{e})))_{w\in
V_n}$, $\ov{e}\in E(v_0,w)$, satisfy the following conditions for $n\geq 1$:

(i) $p^{(n)} \in  C_n^\infty$,

(ii) $F_n^T p^{(n+1)} = p^{(n)}$.


Conversely, if a sequence of vectors $\{ p^{(n)}\}$ from
$\mathbb R^{|V_n|}_+$ satisfies condition (ii), then there
exists a non-atomic  finite Borel $\mathcal R$-invariant measure $\mu$ on
$X_B$ with $p_w^{(n)}=\mu(X_w^{(n)}(\overline e))$ for all $n\geq 1$
and $w\in V_n$.

The $\Rk$-invariant measure $\mu$ is a probability measure if and only if

(iii) $\sum_{w\in V_n}h_w^{(n)}p_w^{(n)}=1$ for $n=1$,

\noindent in which case this equality holds for all $n\ge 1$.
\end{theorem}

{\it Proof}. It follows from (\ref{from})
that if $\mu\in M(\mathcal R)$, then the
sequence $p^{(n)}= (p^{(n)}_w)$  with $p^{(n)}_w = \mu(X^{(n)}_w(\ov{e}))$
satisfies condition (ii). Condition (i) follows from (ii) by the definition of the cones, since
$$
p^{(n)} = F_n^T F_{n+1}^T\cdots F_{n+k}^T p^{(n+k)} \in C_n^{(n+k)}\ \ \mbox{for all}\ k\ge 1.
$$

Conversely, suppose that a sequence of vectors $\{ p^{(n)}\}$
satisfies condition (ii). Define the measure $\mu$ on
$X_w^{(n)}(\overline e),\ w\in V_n,$ to be equal to $p_w^{(n)}$. For
any other clopen set $Y$, we represent $Y$ as a disjoint union of
cylinder sets and define $\mu(Y)$ as the sum of values of $\mu$ on
these cylinder sets. It is routine to check that the measure $\mu$
is well-defined. The definition of $\mu$ yields that it is $\mathcal R$-invariant.
This measure is non-atomic since all the $\Rk$-equivalence classes are infinite by assumption.

The last claim concerning probability measures is immediate since $\{X_w^{(n)}\}_{w\in V_n}$ is a clopen partition
of $X_B$ for any $n\ge 1$. \qed
\medbreak
With every Bratteli diagram $B=(V,E)$ one can associate the
dimension group $ K_0(B)$ \cite{herman_putnam_skau:1992}:
$$
K_0(B)= \varinjlim_n (\Z^{d_n}, F_n) = \Z
\stackrel{F_0}{\to}\Z^{d_1} \stackrel{F_1}{\to} \Z^{d_2}
\stackrel{F_2}{\to} \Z^{d_3} \stackrel{F_3}{\to} \cdots
$$
where $d_n = |V_n|$ and $F_n$ is the incidence matrix. Then $K_0 =
K_0(B)$ is an ordered group whose positive cone $K^+_0$ is naturally
defined by the cones $\Z_+^{d_n}$. Denote by $\underline{1}$ the
ordered unit from $K_0(B)$ corresponding to  $1\in \Z$.

By $S_1(K_0)$, we denote the set of {\it states} on  the dimension
group: $\rho \in S_1(K_0)$ if $\rho$ is a positive homomorphism from
$K_0$ into $\R$ with $\rho(\underline{1}) =1$. The following proposition is known \cite[Theorem 5]{kerov_vershik} (see also 
\cite[p.\,1694]{GJ}),
but we provide a short proof
for the reader's convenience.

\begin{proposition} \label{state} There exists a 1-to-1 correspondence between the sets $M_1(\mathcal R)$ and $S_1(K_0)$.
\end{proposition}
{\it Proof}. We first note that every probability measure on $X_B$ determines uniquely a positive homomorphism on $K_0$.

Conversely, let $\rho : K_0 \to \R$ be a state such that
$\rho(\underline{1}) = 1$. Then there exists a sequence $\rho_i :
\Z^{d_i}\to \R$ of positive homomorphisms such that $\rho_i =
\rho_{i+1}\circ F_i,\ i\geq 0$.
Obviously, $\rho_i(y) = \langle y, \sigma^{(i)} \rangle,\ y
\in \Z^{d_i}$ for some $\sigma^{i} \in \R_+^{d_i}$.
The relation $\rho_i = \rho_{i+1}\circ F_i$ implies that for any
$y\in \Z^{d_i}$,
$$
\rho_{i+1}(F_i y) = \langle F_i y, \sigma^{(i+1)} \rangle = \langle
y, F^T_i\sigma^{(i+1)}\rangle  = \langle y, \sigma^{(i)}\rangle,
$$
hence $F_i^T \sig^{(i+1)}=\sig^{(i)}$ for $i\ge 0$.
By Theorem~\ref{Theorem_measures_general_case}, the sequence $\sig^{(i)}$ determines a measure on $X_B$. This is a probability measure,
because $\rho(\underline{1}) =1$ implies
$$
1 = \rho_0(1) = \rho_1\circ F_0(1) = \rho_1(h^{(1)}) = \langle h^{(1)}, \sig^{(1)} \rangle,
$$
which is the property (iii) of the theorem. \qed

\medskip
\myremark\label{InfiniteMeasure} 1. An
analogue of Theorem
\ref{Theorem_measures_general_case} is valid for the set
$M_\infty(\mathcal R)$ of infinite $\sigma$-finite $\mathcal
R$-invariant measures on the path space $X_B$ of a Bratteli diagram
$B = (V,E)$. Given $\mu \in M_\infty(\mathcal R)$, define   $p^{(n)}
= (p^{(n)}_w)_{w\in V_n}$ where $p^{(n)}_w = \mu(X^{(n)}_w(\ov{e})),\ \ov{e}\in E(v_0,w),\ n \geq
1$. Then at least one of the coordinates of $p^{(n)}$ is
infinite. Relation (\ref{from}) also holds in this
case. More precisely, it shows that if $p^{(n)}_w =\infty,\ w\in
V_n$, then at least one of $p^{(n+1)}_{v_1}, ...,p^{(n+1)}_{v_l}$ is
infinite where $v_1,...,v_l$ are the vertices from $V_{n+1}$ which
are connected with $w$. On the other hand, if  $p^{(n)}_w$ is finite, then all $p^{(n+1)}_{v_1},
...,p^{(n+1)}_{v_l}$ are finite. Conversely, from any sequence of
vectors  $p^{(n)} = (p^{(n)}_w)_{w\in V_n}$ whose coordinates
satisfy the described property one can uniquely restore an infinite
$\mathcal R$-invariant measure.

2. Similarly to Proposition~\ref{state}, the set of infinite $\Rk$-invariant measures corresponds to the set of {\em semi-finite}
states on $K_0$.

3.
After this work was completed, we became aware of the preprints \cite{Fisher1,FFT} where some
related questions are investigated. Fisher \cite{Fisher1} studies minimal non-stationary Bratteli diagrams and
obtains several criteria for unique ergodicity. One of them can be stated as follows, using our
notation:
\begin{quote}
{\em The equivalence relation $\Rk$ on a Bratteli diagram $B$ is uniquely ergodic if and only if the cone $C_n^\infty$ reduces to a single ray for all $n\ge 1$.}
\end{quote}
We note that this is an immediate corollary of Theorem \ref{Theorem_measures_general_case}, and we
do not assume minimality. The proof in \cite{Fisher1} is completely different.
We should  mention that the idea of nested cones was  used by Keane
\cite{keane77} in 1977 to construct a minimal, non-uniquely ergodic
interval exchange transformation. A non-stationary Bratteli-Vershik realization
of Keane's example, as well as several other non-stationary examples of this kind are given
in \cite{FFT}.

%
%
%
%
\section{Non-negative matrices and stationary Bratteli diagrams}\label{stationary}

We first recall some results from the Perron-Frobenius theory of
non-negative matrices. The exposition is based on the papers
\cite{schneider}, \cite{tam-schneider:1994}, and
\cite{tam_schneider:2001}.

Let $F$ be an $N\times N$ matrix with non-negative integer entries
$(f_{i,j})$. Define the directed graph $G(F)$ associated to $F$ whose
vertices are $\{1,...,N\}$ and  there is an arrow from $i$ to  $j$
if and only if $f_{i,j} > 0$. The vertices $i$ and
$j$ are {\it equivalent} if either $ i =j$ or  there is a path in
$G(F)$ from $i$ to  $j$ as well as a path from $j$ to $i$. Let
$\E_i$ denote the corresponding equivalence classes, $i=1,...,m$.
Every class $\E_i$ defines an irreducible submatrix $F_i$ of $F$
obtained by restriction of $F$ to the set of vertices from $\E_i$ (some of $F_i$ may be zero).

Define a partial order on the family of sets $\E_1,..., \E_m$ which
we will identify with $\{1,...,m\}$. For $\alpha,\beta \in
\{1,...,m\}$, we say that a class $\alpha$ {\it has access} to a
class $\beta$, in symbols $\alpha \succeq \beta$, if and only if
either $\alpha=\beta$ or there is a path in $G(F)$ from a vertex
which belongs to $\E_\alpha$ to a vertex which belongs to
$\E_\beta$. We will also say that a vertex $i$ from $G(F)$ {\it is
accessible} from a class $\alpha\in \{1,\ldots,m\}$ if there is a
path in $G(F)$ from a vertex (in fact, from any vertex)  of
$\E_\alpha$  to the vertex $i$. If $\alpha \succeq \beta$ and
$\alpha\neq \beta$, then the notation $\alpha \succ \beta$ is used.
This partial order defines the reduced directed graph $R(F)$ of
$G(F)$ on the set $\{1,...,m\}$ of equivalence classes: by
definition, there is a directed edge in $R(F)$ from $\alpha$ to
$\beta$ if and only if there is a directed edge from a vertex in class $\alpha$ to a vertex in class $\beta$. A vertex $\alpha$ in
$R(F)$ is called {\it final (initial)} if there is no $\beta \in
R(F)$ such that $\beta \prec \alpha$ (respectively $\beta \succ \alpha$).
Slightly different, but equivalent, terminology is used in \cite[4.4]{LM}: $\E_\alpha$ are called there {\em communicating classes},
and initial (final) vertices of the reduced graph are called {\em sources (sinks)} respectively.

One can assume without loss of generality that $\alpha
\succ \beta$ implies that $\alpha > \beta$ (with the usual ordering on integers). Equivalently, the
non-negative matrix $F$ can be transformed by applying
permutation matrices to the Frobenius Normal Form:

\begin{equation}\label{Frobenius Form}
F =\left(
  \begin{array}{ccccccc}
    F_1 & 0 & \cdots & 0 & 0 & \cdots & 0 \\
    0 & F_2 & \cdots & 0 & 0 & \cdots & 0 \\
    \vdots & \vdots & \ddots & \vdots & \vdots & \cdots& \vdots \\
    0 & 0 & \cdots & F_s & 0 & \cdots & 0 \\
    X_{s+1,1} & X_{s+1,2} & \cdots & X_{s+1,s} & F_{s+1} & \cdots & 0 \\
    \vdots & \vdots & \cdots & \vdots & \vdots & \ddots & \vdots \\
    X_{m,1} & X_{m,2} & \cdots & X_{m,s} & X_{m,s+1} & \cdots & F_m \\
  \end{array}
\right)
\end{equation}
The square nonzero matrices $F_\alpha$ standing on the main diagonal are irreducible.  For any fixed $j = s+1,..., m$, at least one of the
matrices $X_{j,k}$ is not zero. Notice that $X_{j,k} \neq 0$ if and only if
there is an edge in $R(F)$ from $\E_j$ to $\E_k$. The fact that
$X_{j,k} = 0$ for all distinct $j,k = 1,...,s$ shows that there are
no edges in $R(F)$ outgoing from the vertices $\{1,...,s\}$, that
is, the vertices $1,...,s$ are final in $R(F)$.
Take the irreducible submatrix $F_\alpha$ corresponding to $\alpha
\in \{1,\ldots,m\}$. Let  $\rho_\alpha= \max\{|\lambda| : \lambda\in
\Spec(F_\alpha)\}$ be the spectral radius of $F_\alpha$. If the
spectrum  $\Spec(F_\alpha)$ contains exactly $h_\alpha$ eigenvalues
$\lambda_1,...,\lambda_{h_\alpha}$ with $|\lambda_i|=\rho_\alpha$,
then $h_\alpha$ is called the index of imprimitivity of $F_\alpha$.
In this case $F^{h_\alpha}_\alpha$ is a primitive matrix. Note that
$F_\alpha$ is primitive if and only if $h_\alpha =1$, see \cite[Section XIII.5]{gant}.

A vertex (class) $\alpha\in \{1,...,m\}$ is called a {\it
distinguished vertex (class)} if $\rho_\alpha > \rho_\beta$ whenever
$\beta \succ  \alpha$. A real number $\lambda$ is called a {\it
distinguished eigenvalue} if there exists a non-negative eigenvector $x$
with $Fx=\lambda x$. Notice that all vertices $\alpha = 1,\ldots, s$ are necessarily distinguished.
The following result extends the well-known Frobenius theorem to the
case of reducible matrices. The proof can be found in
\cite[Proposition 1]{victory}, \cite[Theorem 3.7]{schneider} and
\cite[Theorem 3.3]{tam_schneider:2001}.

\begin{theorem}[Frobenius Theorem] \label{Theorem_FrobeniusTheorem} Let $F$ be an $N\times N$ non-negative matrix with integer entries.

{\bf (a)} A real number $\lambda$ is a distinguished eigenvalue if and
only if there exists a distinguished class $\alpha$ in $R(F)$ such
that $\rho_\alpha=\lambda$.

{\bf (b)} If $\alpha$ is a distinguished class in $R(F)$, then there
exists a unique (up to scaling) non-negative eigenvector
$\xi_\alpha=(x_1,\ldots,x_N)^T$ corresponding to $\rho_\alpha$ having the property
that  $x_i>0$ if and only if the vertex $i$ has access to $\alpha$. 
\end{theorem}

Note that in part (b) the uniqueness refers to eigenvectors with the given property; there may be other non-negative eigenvectors corresponding to 
$\rho_\alpha$ if there is another distinguished class with the same spectral radius (these classes will be necessarily non-accessible to each other).

We will call $\xi_\alpha$ from
Theorem~\ref{Theorem_FrobeniusTheorem} the {\em distinguished
eigenvector} corresponding to $\alpha$.

For a non-negative $N\times N$ matrix $A$, define
$$
core(A)=\bigcap_{k\geq 1}A^k(\mathbb R^N_+).
$$
For a non-negative matrix $A$ and $k\in \N$, denote by $C(A,k)$ the
cone generated by the distinguished eigenvectors of $A^k$. Let
$\Lambda\subset \{1,...,m\}$ be the set of all irreducible
components of $A$ with  positive spectral radii. For every
$\alpha\in\Lambda$ denote by $h_\alpha$  the index of imprimitivity
of  the irreducible component $A_\alpha$.

The following theorem (see  \cite[Theorem 4.2]{tam-schneider:1994})
describes $core(A)$ for non-negative matrices. This result is of
crucial importance for our study of invariant measures.

\begin{theorem}\label{Theorem_Core} Let $A$ be a non-negative  $N\times N$ matrix with positive spectral radius. Then

{\bf (a)} $core(A)$ is a simplicial cone with exactly
$\sum_{\alpha\in\Lambda}{h_\alpha}$ extreme rays.

{\bf (b)} $core(A)=C(A,q)$ where $q$ is the least common multiple of all
$h_\alpha$, $\alpha\in \Lambda$. In particular, if all the irreducible components of $A$ are primitive, then $core(A) = C(A,1)$.
\end{theorem}
\medbreak

\myremark\label{RemarkMatrixStructure}  If $B=(V,E)$ is a stationary
Bratteli diagram, then we have {\it two} non-negative integer
matrices associated to $B$: the incidence matrix $F$ and its
transpose matrix $A$. The reduced graphs $R(F)$ and $R(A)$ have the
same sets of vertices but the opposite direction of edges. This
means that if $\alpha \succeq \beta$ for $R(F)$, then $\alpha
\preceq\beta$ in $R(A)$. Saying that
 $\alpha$ has access to $\beta$, we need to point out the graph, $R(F)$ or $R(A)$, in which these vertices are considered.    It follows that  the reduced graphs have different sets of distinguished vertices.  More precisely, if $F$ is represented in Frobenius form (\ref{Frobenius Form}), then
\begin{equation}\label{Frobenius Form for A}
A =\left(
  \begin{array}{ccccccc}
    A_1 & 0 & \cdots & 0 & Y_{1,s+1} & \cdots & Y_{1,m} \\
    0 & A_2 & \cdots & 0 & Y_{2,s+1} & \cdots & Y_{2,m} \\
    \vdots & \vdots & \ddots & \vdots & \vdots & \cdots& \vdots \\
    0 & 0 & \cdots & A_s & Y_{s,s+1} & \cdots & Y_{s,m} \\
    0 & 0 & \cdots & 0 & A_{s+1} & \cdots & Y_{s+1,m} \\
    \vdots & \vdots & \cdots & \vdots & \vdots & \ddots & \vdots \\
    0 & 0 & \cdots & 0 & 0 & \cdots & A_m \\
  \end{array}
\right)
\end{equation}
where $A_i$ is the matrix transpose to $F_i$ and $Y_{j,i}$ is
transpose to $X_{i,j}$. The vertices $\{1,...,s\}$ in the graph
$R(F)$ are final but the same vertices in $R(A)$ are initial. Then
we obtain, in particular, that  $\{1,...,s\}$ are distinguished
vertices in $R(A)$.

Now suppose that $B = (V,E)$ is a stationary Bratteli diagram with the incidence matrix $F$.
Fix $d\in \N$ and consider the Bratteli diagram $B_d $ which is
obtained by telescoping $B$ with respect to the levels $V_{nd+1},\ n
= 0,1,\ldots$, see \cite{herman_putnam_skau:1992}. Then $B_d$ is again a stationary diagram, whose
incidence matrix is $F^d$. There is an obvious way to identify the path spaces $X_{B_d}$ and $X_B$,
which preserves the tail equivalence relation. Therefore, we can naturally identify the invariant measures
for these diagrams. Thus,
without loss of generality, we can telescope the diagram $B$ and
regroup the vertices in such a way that the matrix $F$ will have the
following property:
\begin{eqnarray} \label{prop-I}
& F & \mbox{has the form (\ref{Frobenius Form}) where every nonzero matrix $F_i$}\\[-1ex]
 &  & \mbox{on the main diagonal  is primitive.} \nonumber
\end{eqnarray}
We can telescope the diagram $B$ further to make sure that
\begin{eqnarray} \label{prop-I'}
& F & \mbox{has the form (\ref{Frobenius Form}) where every nonzero matrix $F_i$}\\[-1ex]
& & \mbox{on the main diagonal  is strictly positive.} \nonumber
\end{eqnarray}
However, this is not always convenient, since it may lead to
matrices with large entries.
We record some properties of the matrix $F$ in the next lemma.

\begin{lemma} \label{lem-prop} Suppose that $B$ is a stationary Bratteli diagram such that the 
tail equivalence relation $\Rk$ is aperiodic and
the incidence matrix $F$ satisfies (\ref{prop-I}). Then

{\bf (a)} $F_i \ne 0$ for $1\le i \le s$.

{\bf (b)} If $F_i$ is a non-zero $1\times 1$ matrix, then its entry is greater than one, for $1\le i \le s$.

{\bf (c)} $core(A)$, where $A=F^T$, is the simplicial cone generated by the distinguished eigenvectors of $A$.
All distinguished eigenvalues of $A$ are greater than one.
\end{lemma}

Note that we do not exclude that some of the matrices $F_i$, for $s < i \le m$, are $1\times 1$ of the form $[0]$
or $[1]$.

\medskip

{\em Proof.} (a) By Definition \ref{Definition_Bratteli_Diagram}, there are edges leading into every vertex.
Since the vertices in the classes $1\le i \le s$ are initial in the graph $R(A)$, there have to be some edges from class $i$ to itself.
The claim follows.

(b) If a class $i$, for $1\le i \le s$, consists of only one vertex with only one loop edge, then it defines an infinite path in the
diagram such that its $\Rk$-equivalence class consists of a single element (again, because it is an initial vertex in $R(A)$), which
contradicts our assumption that it is infinite.

(c) The first statement is contained in Theorem~\ref{Theorem_Core}(b). Now let $\alpha$ be a distinguished vertex of $A$. By definition,
the distinguished eigenvalue $\lam_\alpha$ is greater than all PF eigenvalues of classes which have access to $\alpha$ in
$R(A)$. There is always such a class which is an initial vertex of $R(A)$, and its PF eigenvalue is greater than one by parts (a) and (b)
of this  lemma. It follows that $\lam_\alpha>1$. \qed

\medskip

Denote by $\Lambda$ the set of those vertices $\alpha$ in $R(A)$
for which $F_\alpha \neq 0$, and let $A_\alpha=F_\alpha^T$.
For $\alpha\in \Lambda$, denote by $B_\alpha$ the stationary
subdiagram of $B$ consisting of vertices which belong to the class
$\E_\alpha$ and those edges which connect them. 
Condition (\ref{prop-I}) means that the subdiagram $B_\alpha$ is
simple. 

Let $Y_\alpha$ be the path space of the Bratteli diagram $B_\alpha,\
\alpha \in \Lambda$.  Define $X_\alpha = \mathcal
R(Y_\alpha)$, that is, a path $x\in X_B$ belongs to $X_\alpha$ if it
is $\mathcal R$-equivalent to a path $y\in Y_\alpha$. It is clear that 
$\{X_\alpha : \alpha \in \Lambda\}$ is a partition of $X_B$.
The following lemma describes the orbit closures for the equivalence relation and the minimal components.

\begin{lemma} \label{lem-top}
Under the assumption (\ref{prop-I}) we have
$$
\forall\,\alpha\in \Lam,\ \forall\,x\in X_\alpha,\ \ov{\Rk(x)} = \bigcup_{\beta\in \Ik_\alpha} X_\beta,\ \ \mbox{where}\ \Ik_\alpha = \{\beta\in \Lam:\ 
\beta \succeq \alpha\}.
$$
Thus, the minimal components of $X_B$ for the tail equivalence relation are exactly $X_\alpha,\ \alpha =1,\ldots,s$, that is, those
$X_\alpha$ which corresponding to the
initial vertices of $R(A)$.
\end{lemma}

{\em Proof} is immediate from the structure of the diagram and the definitions. \qed

\medskip

Next we obtain  necessary and sufficient conditions under which
a measure on $X_B$ is $\mathcal R$-invariant. We use the
notation of Section \ref{section2}. Recall that we can assume property (\ref{prop-I}) without loss of generality.

\begin{theorem}\label{TheoremStationaryMeasures}
 Suppose that $B$ is a stationary Bratteli diagram such that the
tail equivalence relation $\Rk$ is aperiodic and
the incidence matrix $F$ satisfies (\ref{prop-I}).
Let $\mu$ be a finite Borel $\mathcal R$-invariant measure on $X_B$.
Set $p^{(n)}=(\mu(X_w^{(n)}(\overline e)))_{w\in V_n}$ where  $\overline e\in E(v_0,w)$. Then for the  matrix $A = F^T$ associated to $B$ the following properties hold:

(i) $p^{(n)}=A p^{(n+1)}$ for every $n\geq 1$;

(ii) $ p^{(n)}\in core(A),\ n\geq 1$.


Conversely, if a sequence of vectors $\{ p^{(n)}\}$ from
$\mathbb R^{N}_+$ satisfies condition (ii), then there
exists a finite Borel $\mathcal R$-invariant measure $\mu$ on
$X_B$ with $p_w^{(n)}=\mu(X_w^{(n)}(\overline e))$ for all $n\geq 1$
and $w\in V_n$.

The $\Rk$-invariant measure $\mu$ is a probability measure if and only if

(iii) $\sum_{w\in V_n}h_w^{(n)}p_w^{(n)}=1$ for $n=1$,

\noindent in which case this equality holds for all $n\ge 1$.
\end{theorem}

{\it Proof.} This is just a special case of Theorem~\ref{Theorem_measures_general_case}. We only need to note that
$C_n^\infty = core(A)$ when $B$ is a stationary diagram. \qed

\medskip
The next lemma together with Theorem \ref{TheoremStationaryMeasures}
shows that each vector  $ p^{(1)}\in core(A)$ uniquely defines a finite
$\Rk$-invariant probability measure $\mu$ on $X_B$. We note that
this lemma follows implicitly from \cite[Theorem
2.2]{tam-schneider:1994}.

\begin{lemma}\label{Lemma_Uniqueness} Let $\{p^{(n)}\}_{n\geq 1}$ and $\{ q^{(n)}\}_{n\geq 1}$ be two sequences of vectors in $\mathbb R^N$ such that $ p^{(n)}=A p^{(n+1)}$ and $ q^{(n)}=A q^{(n+1)}$ for all $n\geq 1$. If $ p^{(1)}= q^{(1)}$, then $ p^{(n)}=  q^{(n)}$ for every $n\geq 1$.
\end{lemma}
{\it Proof.} Suppose this is not true. Take the first integer $n_0\geq 1$
with $ p^{(n_0+1)}\neq  q^{(n_0+1)}$. Clearly,
$ p^{(n)}\neq q^{(n)}$ for every $n> n_0 +1$. For
each $j\geq 1$, set $ x^{(j)}=
p^{(n_0+j)}- q^{(n_0+j)}\neq 0$. It follows that
$A^j x^{(j)}=0$ whereas $A^{j-1} x^{(j)}=
x^{(1)}\neq 0$. This implies that the family $\{
x^{(1)},\ldots,  x^{(j)}\}$ is linearly independent for any
$j\geq 1$, which is impossible.  \hfill$\square$ \medbreak

Consider the cone $core(A)$. Denote by $\xi_1,\ldots,\xi_k$ the
extreme vectors of $core(A)$.  We  normalize each vector $\xi_i$ so
that $\sum_{w\in V_1}h_w^{(1)}(\xi_i)_w=1$.  Then each vector of
$D=\{ x\in core(A) : \sum_{w\in V_1}h_w^{(1)}x_w=1\}$ is a
convex combination of the vectors $\xi_1,\ldots,\xi_k$.  The next
theorem is one of the main result of this paper  which completely
describes the simplex of $\mathcal R$-invariant probability measures
of a  stationary Bratteli diagram.

\begin{theorem}\label{Theorem_MeasuresErgodicMeasures}
 Suppose that $B$ is a stationary Bratteli diagram such that the
tail equivalence relation $\Rk$ is aperiodic and
the incidence matrix $F$ satisfies (\ref{prop-I}).
Then there is a one-to-one correspondence between vectors $p^{(1)}\in D$ and $\mathcal R$-invariant probability
measures on $X_B$. This correspondence is given by the rule $\mu \leftrightarrow
p^{(1)}= (\mu(X_w^{(1)})/h_w^{(1)})_{w\in V_1}$. Furthermore, ergodic measures correspond to the extreme vectors $\{\xi_1,\ldots,\xi_k\}$. In particular, there exist exactly $k$ ergodic measures.
\end{theorem}
{\it Proof.}  By Lemma~\ref{lem-prop}(c), we have that every extreme vector
$\xi_i$ is an eigenvector of $A$ for some distinguished  eigenvalue
$\lambda_i>1$, $i=1,\ldots,k$.

Take a vector $p^{(1)}\in D$.
By the definition of $D$ we can find a sequence $p^{(n)}\in core(A)$ such that $p^{(n)} = A p^{(n+1)}$
for every $n\ge 1$. Thus, this sequence satisfies conditions (i) and (ii) of
Theorem~\ref{TheoremStationaryMeasures} and hence
defines an
$\mathcal R$-invariant finite Borel measure $\mu$ on $X_B$. Condition (iii) also holds, whence $\mu$ is a probability measure.
It
immediately follows from Lemma \ref{Lemma_Uniqueness} and Theorem
\ref{TheoremStationaryMeasures} that there exists only one measure
$\mu$ with   $p^{(1)}=(\mu(X_w^{(1)})/h_w^{(1)})_{w\in V_1}$.

Observe that the correspondence $\mu \leftrightarrow p^{(1)}$ is affine linear, therefore, the simplex
of $\Rk$-invariant probability measures is affine-homeomorphic to $D$. It is well-known that ergodic
invariant measures are precisely the extreme points of this simplex (see \cite[Theorem 6.10]{walters:book}),
which yields the last claim of the theorem. \qed

\medskip

\myremark\label{measureCylinderSets} 1. Observe that if $\mu_\alpha$
is the ergodic measure corresponding to the distinguished
eigenvector $\xi_\alpha = (\xi_\alpha(1),..., \xi_\alpha(N))^T$ then
there is a simple formula for computing the measure $\mu_\alpha$ of
cylinder sets. Let $X_v^{(n)}(\ove)$ be the cylinder set defined by
the finite path $\ove = (e_1,...,e_n)$ with $r(e_n) = v$. Then
\begin{equation}\label{measure_cylinder_sets}
\mu_\alpha(X_v^{(n)}(\ove)) =
\frac{\xi_\alpha(v)}{\lambda_\alpha^{n-1}}.
\end{equation}

2. If $A$ is a primitive matrix, then $core(A)$ is exactly
the ray generated by the Perron-Frobenius eigenvector. By Theorem
\ref{Theorem_MeasuresErgodicMeasures} we get another proof of the well-known fact that
the Vershik map on a stationary Bratteli diagram with a primitive incidence matrix is uniquely ergodic.

3. Handelman \cite{Handelman} studied the dimension group of reducible Markov chains. His \cite[Theorem I.3]{Handelman} is similar to
our Theorem~\ref{Theorem_MeasuresErgodicMeasures}, in the setting of states on the dimension group, although the proof is
completely different. However, his characterization of non-negative eigenvectors \cite[Theorem I.1]{Handelman} is incorrect, except in
the 2-component case. The main part of \cite{Handelman} is devoted to the 2-component case and the description of the dimension group as
as an extension, in terms of the dimension groups of the irreducible components.

%
%
\section{Finite and infinite invariant measures}
%

Here we obtain some additional properties of ergodic $\Rk$-invariant probability measures
described in the previous section and then characterize infinite ($\sigma$-finite) non-atomic invariant
measures.

Let $B$ be a stationary Bratteli diagram. Recall that we can assume the property (\ref{prop-I}) without loss of generality, and this
will be a standing assumption throughout this section.
In Theorem
\ref{Theorem_MeasuresErgodicMeasures} we obtained a complete
description of the set of ergodic probability $\mathcal R$-invariant
measures on the path space $X_B$. Let
$\alpha$ be a distinguished vertex of the reduced graph $R(A)$
with the vertex set $\{1,...,m\}$, and let $\lambda_\alpha =
\rho(A_\alpha)$ be the Perron-Frobenius eigenvalue of $A_\alpha$. Recall that
the corresponding distinguished eigenvector $\xi_\alpha =
(\xi_1,...,\xi_N)^T$ of the matrix $A$ has the property $\xi_i >0$ if and
only if the vertex $i$ has access to $\alpha$.  Recall also that
$\lambda_\alpha > \rho(A_\beta)$ for every class $\beta$ which has
access to $\alpha$ in $R(A)$. Below we write $\beta \succeq \alpha$
if $\beta$ has access to $\alpha$ in the graph $R(A)$.

It follows from \cite[Theorem 9.4]{schneider} that
\begin{equation} \label{eq-asymp}
(A^n)_{i,j} \sim \lam_\alpha^n,\ \ n\to \infty,\ \ \mbox{for}\ i\in
\Ek_\beta,\ j\in \Ek_\alpha,\ \mbox{with}\ \beta \succeq \alpha.
\end{equation}
Here $\sim$ means that the ratio tends to a positive constant. On
the other hand,
\begin{equation} \label{eq-asymp2}
(A^n)_{i,j} = o(\lam_\alpha^n),\ \ n\to \infty,\ \ \mbox{for}\ j\in
\Ek_\beta,\ \mbox{and any $i$,\ \ with}\ \beta\succ \alpha.
\end{equation}
(There is precise asymptotics, depending on $i$, in the latter case as well, but we do
not need it.)

%
%

Recall the notation introduced in Section 3 after Lemma~\ref{lem-prop}: the set
$\Lam=\{\beta:\ A_\beta \ne 0\}$, the simple subdiagram $B_\beta$ corresponding to $\beta\in \Lam$, and the partition $\{X_\beta:\ \beta\in \Lam\}$ of $X_B$.
Let $\ov{e} = (e_1,\ldots,e_m)$ be a finite path in $B$ from $v_0$
to the level $m$; recall the notation $[\ov{e}] = X_v^{(m)}(\ov{e})$, with $v =
r(e_m)$. For a finite path $\ov{\om} =
(\om_1,\ldots,\om_m)$ in $B$ (not necessarily starting from $v_0$)
we write $s(\ovom) = s(\om_1)$ and $r(\ovom) = r(\om_m)$.

Fix the ergodic $\mathcal R$-invariant probability measure
$\mu_\alpha$ corresponding to a distinguished vertex $\alpha$ of the
reduced graph $R(A)$.

\begin{lemma} \label{lem-support} For a distinguished vertex $\alpha$, the measure $\mu_\alpha$ is supported on $X_\alpha$.
\end{lemma}
{\it Proof.} Let $\xi_\alpha = (\xi_\alpha(1),...,\xi_\alpha(N))^T$
be the distinguished eigenvector corresponding to $\alpha$. To prove
the lemma, it is enough to show that $\mu_\alpha(X_\beta) =0$ for $\beta\in \Lam$, 
$\beta\ne \alpha$. If $\beta$ does not have access to $\alpha$ in
$R(A)$, then this is immediate, since for every finite path $\ove$
with $r(\ove)\in V(B_\beta)$ we have $\mu_\alpha([\ove])=0$ (see
Theorems \ref{Theorem_FrobeniusTheorem},
\ref{Theorem_MeasuresErgodicMeasures} and Remark
\ref{measureCylinderSets}). Now suppose that $\beta$ has access to
$\alpha$ in $R(A)$. We can write $X_\beta = \bigcup_{\ell\ge
1}X_\beta^{(\ell)}$, where $X_\beta^{(\ell)}$ is the set of $x
=(x_n)\in X_\beta$ such that $x_n\in E(B_\beta)$ for $n\ge \ell$,
and prove that $\mu_\alpha(X_\beta^{(\ell)})=0$ for all $\ell$.
Recall that for every finite path $\ove$ of length $n$ we have
$\mu_\alpha([\ove]) = \xi_\alpha(v) \lam_\alpha^{-n+1}$ where $v
=r(\ove)$ by (\ref{measure_cylinder_sets}). The number of paths of
length $n$ which terminate in $V(B_\beta)$ equals
$$
\sum_{j\in \Ek_\beta} h^{(n)}_j =\sum_{j\in \Ek_\beta}\sum_{i=1}^N
((A^T)^{n-1})_{j,i}h^{(1)}_i = \sum_{j\in \Ek_\beta}\sum_{i=1}^N
(A^{n-1})_{i,j}h^{(1)}_i,
$$
which is $o(\lam_\alpha^{n-1})$ in view of (\ref{eq-asymp2}). Notice
also that for any $n\geq \ell$ we have that
$$\mu_\alpha(X_\beta^{(\ell)})\leq \sum_{j\in
\Ek_\beta}\frac{\xi_\alpha(j)h_j^{(n)}}{\lambda_\alpha^{n-1}}\leq
\frac{1}{\lambda_\alpha^{n-1}} \sum_{j\in \Ek_\beta} h^{(n)}_j.$$
Since we can choose $n$ arbitrarily large, it follows that
$\mu_\alpha(X_\beta^{(\ell)})=0$.  \hfill$\square$
\medskip

If $A_\alpha \neq 0$, then there exists a unique $\mathcal
R_\alpha$-invariant probability measure $\nu_\alpha$ on the path
space $Y_\alpha$ of $B_\alpha$ where $\mathcal R_\alpha = \mathcal R
\cap (Y_\alpha \times Y_\alpha)$. We can naturally extend the
measure $\nu_\alpha$ to the space  $X_\alpha$ and produce there a
measure $\widetilde \nu_\alpha$ which is $\mathcal R$-invariant. In
fact, $X_\alpha \setminus Y_\alpha$ is a disjoint union of cylinder
sets $[\ove]$ corresponding to paths $\ove = (e_1,\ldots,e_m)$ for
some $m\ge 1$, such that $r(e_m) \in V(B_\alpha)$, but $s(e_m)
\not\in V(B_\alpha)$. For each such cylinder set the measure
$\nutil|_{[\ove]}$ is defined to be a copy of $\nutil|_{[\ove']} =
\nu_{\alpha}|_{[\ove']}$ for a path $\ove' = (e'_1,\ldots, e'_m)\in
B_\alpha$ with $r(e_m)=r(e'_m)$. Observe that if we equip the
Bratteli diagram $B$ with an order, then it defines an order on
$B_\alpha$. Let $\varphi_B$ and  $\varphi_\alpha$ be the Vershik
maps defined on $X_B$ and $X_{B_\alpha}= Y_\alpha$ (with the orbits
of maximal and minimal paths removed). Recall that $B_\alpha$ is a
simple diagram ($A_\alpha$ is primitive), hence
$(Y_\alpha,\varphi_\alpha)$ is uniquely ergodic. Therefore, the
measure $\nutil_\alpha$ is an ergodic (possibly infinite) measure
for the induced transformation $\varphi_B$ (see e.g. \cite[Exercise
1, p.\,56]{petersen:book}).

In the next lemma we describe infinite ergodic $\mathcal
R$-invariant measures on the path space of a stationary diagram and
clarify the relation between measures $\mu_\alpha$ and
$\widetilde\nu_\alpha$.

\begin{lemma} \label{lem-infmeas}
 Suppose that $B$ is a stationary Bratteli diagram such that the
tail equivalence relation $\Rk$ is aperiodic and
the incidence matrix $F$ satisfies (\ref{prop-I}).
Suppose $\alpha$ is
a vertex in the reduced graph $R(A)$. If $\alpha$ is a distinguished
vertex, then $\nutil_\alpha = c_\alpha \mu_\alpha$ for some
$c_\alpha>0$. If $\alpha$ is not a distinguished vertex, then
$\widetilde\nu_\alpha $ is an infinite ergodic $\mathcal
R$-invariant measure. The measure $\nutil_\alpha$ is non-atomic, unless $A_\alpha$ is the $1\times 1$ matrix $[1]$.
Conversely, every infinite ergodic invariant
measure which is positive and finite on at least one open set (depending on the measure)
equals $c\nutil_\alpha$ for some $c>0$ and some non-distinguished vertex $\alpha$.
\end{lemma}
{\it Proof}. If $\alpha$ is a distinguished vertex, then the measure
$\mu_\alpha$ is $\varphi_B$-invariant and positive on the cylinders
of $Y_\alpha$. Then $\mu_\alpha|_{Y_\alpha}$ is positive and
invariant for the first return map $\varphi_\alpha$ on $Y_\alpha$,
which is uniquely ergodic. It follows that $\nu_\alpha = c_\alpha
\mu_\alpha|_{Y_\alpha}$, and hence $\nutil_\alpha =
c_\alpha\mu_\alpha$.

If $\alpha$ is not a distinguished vertex, then $\nutil_\alpha$
cannot be finite, since this would contradict
Theorem~\ref{Theorem_MeasuresErgodicMeasures}. If $\nu_\alpha$ is non-atomic, then its extension $\nutil_\alpha$ is
non-atomic. This holds for any non-zero component $A_\alpha$, except when $A_\alpha = [1]$. In the latter case, $Y_\alpha$ is a
singleton, hence $\nu_\alpha$ is a point mass, and its extension $\nutil_\alpha$ is a pure discrete $\sigma$-finite measure.

It remains to verify the last statement of the theorem. Let $\mu$ be
an infinite ergodic $\varphi_B$-invariant measure, which is positive
and finite on an open set. Since cylinder sets generate the
topology, we can find $\ove=(e_1,\ldots,e_m)$ such that $0<
\mu([\ove]) < \infty$. Clearly $r(e_m)\in B_\alpha$ for some
$\alpha$. Without loss of generality, we can assume that $\alpha$ is
the largest index which appears this way. This means that if $\ove$
can be prolonged to a path $\ove'$ with a terminal vertex in another
subdiagram $B_\beta$, $\beta\ne \alpha$, then $\mu([\ove'])=0$. Note
that $\mu|_{Y_\alpha}$ is a finite positive
$\varphi_\alpha$-invariant measure, hence $\nu_\alpha =
c_\alpha\mu|_{Y_\alpha}$ for some $c_\alpha>0$, because $B_\alpha$
is uniquely ergodic, being a simple diagram. Then necessarily $\nutil = c_\alpha\mu$, since two ergodic
(finite or infinite) measures that agree on a set of positive measure are equal.
\hfill$\square$
\medskip

The first part of Lemma \ref{lem-infmeas} can also be proved in a different way. To
show that $\nutil_\alpha$ is finite (and therefore proportional to
$\mu_\alpha$) when $\alpha$ is a distinguished vertex, we can
compute the $\nutil_\alpha$ measure of the set $X_\alpha(n):= \{x
=(e_n)\in X_\alpha : r(e_m)\in V(B_\alpha), m\geq n\}$ for any $n$.
Then
\begin{equation}\label{finite-measure}
\nutil_\alpha(X_\alpha(n)) = \sum_{v\in V_n(B_\alpha)}
h_v^{(n)}\nu_\alpha([\ove_v])
\end{equation}
where $\ove_v$ is a finite path in $B_\alpha$ connecting $v_0$ and
$v$. Since $\nu_\alpha([\ove_v]) = c \lambda_\alpha^{n-1}, c >0$, we
can apply (\ref{heights}) and (\ref{eq-asymp}) to deduce that
$\nutil_\alpha(X_\alpha(n))$ is finite and independent of $n$. If
the vertex $\alpha$ is not distinguished, then $\lambda_\alpha \leq
\lambda_\beta$ for some vertex $\beta$ which has access to $\alpha$
in $R(A)$. Then $h_v^{(n)}$ will grow as $\lambda_\beta^n$ and
$\nutil_\alpha(X_\alpha(n))$ in (\ref{finite-measure}) tends to
infinity as $n\to\infty$.

Thus, we obtained the following result on infinite $\Rk$-invariant measures for stationary Bratteli diagrams.

\begin{theorem} \label{th-infmeas}
 Suppose that $B$ is a stationary Bratteli diagram such that the
tail equivalence relation $\Rk$ is aperiodic and
the incidence matrix $F$ satisfies (\ref{prop-I}).
Then the set of ergodic infinite ($\sig$-finite) invariant measures, which are positive
and finite on at least one open set (depending on the measure), modulo a
constant multiple, is in 1-to-1 correspondence with the set of non-distinguished vertices of the reduced graph $R(A)$, for $A=F^T$.
\end{theorem}

%
%

\section{Applications and Examples}

\subsection{Orbit equivalence.}
 Recall that two topological dynamical systems
$(X_1,T_1)$ and $(X_2,T_2)$ are orbit equivalent if there is a homeomorphism $f:\,X_1\to X_2$ which sends $T_1$-orbits
into $T_2$ orbits.
Giordano, Putnam, and Skau proved in
\cite{giordano_putnam_skau:1995} (among other things) the following result: Let
$(X_1,T_1)$ and $(X_2,T_2)$ be uniquely ergodic minimal
homeomorphisms of Cantor sets and let $\mu_1$ and $\mu_2$ be $T_1$-
and $T_2$-invariant probability measures, respectively. Then $(X_1,T_1)$ and
$(X_2,T_2)$ are orbit equivalent if and only if $\{\mu_1(E):\ E\
\mbox {clopen in } X_1\} = \{\mu_2(E):\ E\ \mbox {clopen in }
X_2\}$. We are going to show that this statement is not valid for
two non-minimal aperiodic uniquely ergodic homeomorphisms.

Let $B_1$ and $B_2$ be two stationary Bratteli diagrams constructed
by the incidence matrices $F_1$ and $F_2$ where
$$
F_1 = \left(
  \begin{array}{cc}
    2 & 0 \\
    1 & 2 \\
  \end{array}
\right), \ \ \ \ F_2 = \left(
                        \begin{array}{ccc}
                          2 & 1 & 0 \\
                          0 & 2 & 0 \\
                          0 & 1 & 2 \\
                        \end{array}
                      \right).
$$
\begin{center}
\begin{tabular*}{0.75\textwidth}%
     {@{\extracolsep{\fill}}cc}
\unitlength = 0.3 cm

 \begin{graph}(8,15)
\graphnodesize{0.4}
 \roundnode{V0}(4.0,14)
\put(0.7,12.5){$C_1$}
\put(6,12.5){$D_1$}
\roundnode{V11}(0.0,9) \roundnode{V12}(8.0,9)
\roundnode{V21}(0.0,5) \roundnode{V22}(8.0,5)
\roundnode{V31}(0.0,1) \roundnode{V32}(8.0,1)
\edge{V11}{V0} \edge{V12}{V0}
\bow{V21}{V11}{0.06} \bow{V21}{V11}{-0.06} \edge{V22}{V11}
\bow{V22}{V12}{0.06} \bow{V22}{V12}{-0.06}
\bow{V31}{V21}{0.06} \bow{V31}{V21}{-0.06} \edge{V32}{V21}
\bow{V32}{V22}{0.06} \bow{V32}{V22}{-0.06}
\freetext(4,-.5){$.\ .\ .\ .\ .\ .\ .\ .\ .\ .\ $}%
 \freetext(4,-2.5){{Diagram $B_1$}}
 \freetext(17,-3.0){{Fig. 2}}
\end{graph}
            &
\unitlength = 0.3 cm
 \begin{graph}(11,15) \graphnodesize{0.4}
 \roundnode{V0}(5.5,14)
 \put(4,11){$C_2$}
 \put(0,11){$D_2$}
 \put(9,11){$D_2$}
\roundnode{V11}(0.0,9) \roundnode{V12}(5.5,9)
\roundnode{V13}(11.0,9)
\roundnode{V21}(0.0,5) \roundnode{V22}(5.5,5)
\roundnode{V23}(11.0,5)
\roundnode{V31}(0.0,1) \roundnode{V32}(5.5,1)
\roundnode{V33}(11.0,1)
\edge{V11}{V0} \edge{V12}{V0} \edge{V13}{V0}
\bow{V21}{V11}{0.06} \bow{V21}{V11}{-0.06} \edge{V21}{V12}
\bow{V22}{V12}{0.06} \bow{V22}{V12}{-0.06} \edge{V23}{V12}
\bow{V23}{V13}{0.06} \bow{V23}{V13}{-0.06}
\bow{V31}{V21}{0.06} \bow{V31}{V21}{-0.06} \edge{V31}{V22}
\bow{V32}{V22}{0.06} \bow{V32}{V22}{-0.06} \edge{V33}{V22}
\bow{V33}{V23}{0.06} \bow{V33}{V23}{-0.06}
\freetext(5.5,-.5){$.\ .\ .\ .\ .\ .\ .\ .\ .\ .\ .\ .\ .\ $}%
 \freetext(5.5,-2.5){{Diagram $B_2$}}
\end{graph}
\\
\vspace{1cm}
\end{tabular*}
\end{center}

It is not hard to see that these diagrams admit Vershik maps
$\varphi_1$ and $\varphi_2$ acting on the path spaces $X_{B_1}$ and
$X_{B_2}$, respectively. Then the dynamical systems
$(X_{B_1}, \varphi_1)$ and $(X_{B_2}, \varphi_2)$ are aperiodic and each has a unique
minimal component: $C_1 \subset X_{B_1}$, corresponding to the left
part of diagram $B_1$ and $C_2\subset X_{B_2}$, corresponding to
the central part of diagram $B_2$. It follows from our results in Section 3 that the systems
$(X_1,T_1)$ and $(X_2,T_2)$ are uniquely ergodic.
Let $\mu_1$ and $\mu_2$ be the
unique ergodic invariant probability measures, which are supported on
$C_1$ and $C_2$, respectively. Notice that $(C_1, \varphi_1)$ and $(C_2, \varphi_2)$ are
identical; in fact, this is the 2-odometer. Therefore the values
of these measures on clopen subsets in $X_{B_1}$ and $X_{B_2}$ are
the same.

\begin{proposition}  The Vershik maps $\varphi_1$ and $\varphi_2$ are not orbit equivalent.
\end{proposition}

{\it Proof}. Suppose that the systems $(X_{B_1}, \varphi_1)$ and
$(X_{B_2}, \varphi_2)$ are orbit equivalent. Notice that a
homeomorphism $f$ implementing orbit equivalence maps $C_1$ onto
$C_2$. Therefore, $(X_{B_1} \setminus C_1, \varphi_1)$ is orbit
equivalent to $(X_{B_2} \setminus C_2, \varphi_2)$ via $f$. Let $D_1$
and $D_2$ be clopen subsets of $X_{B_1}$ and $X_{B_2}$  defined as
shown on Figure 2. It is obvious that $D_1$ and $D_2$ are complete
sections for $(X_{B_1} \setminus C_1, \varphi_1)$ and $(X_{B_2} \setminus C_2, \varphi_2)$ respectively. Let $\psi_1 = (\varphi_1)_{D_1}$ and
$\psi_2 = (\varphi_2)_{D_2}$ be the induced homeomorphisms defined
on $D_1$ and $D_2$. It is straightforward to check that orbit
equivalence of  $(X_{B_1} \setminus C_1, \varphi_1)$ and $(X_{B_2} \setminus
C_2, \varphi_2)$ implies orbit equivalence of $(D_1, \psi_1)$ and
$(D_2,\psi_2)$. But the latter is impossible because $(D_1, \psi_1)$
is a uniquely ergodic system and $(D_2,\psi_2)$ has two ergodic
invariant probability measures. The proposition is proved.
\hfill$\square$
\\

We note that one can use another argument to prove the proposition. It follows from Theorem \ref{th-infmeas} that the diagrams $B_1$ and $B_2$ have different numbers of (essentially distinct) infinite invariant measures.

\medskip

We can apply our results to a question on orbit equivalence in Borel dynamics.

\begin{corollary}\label{Borel Orbit Equiv}
Let $B_1$ and $B_2$ be stationary Bratteli diagrams with incidence matrices $F_1$ and $F_2$ respectively. The tail equivalence
relations $\Rk_1$ and $\Rk_2$ are Borel isomorphic if and only if the matrices $A_1 = F_1^T$ and $A_2=F_2^T$
have the same number of distinguished eigenvalues.

In particular, if
$\omega$ and $\omega'$ are two orderings on a stationary Bratteli diagram $B$ such that Borel-Vershik automorphisms $f_\omega$ and $f_{\omega'}$ of the path space $X_B$ exist, $f_\omega$ and $f_{\omega'}$ are (Borel) orbit equivalent.
\end{corollary}

{\it Proof.} By \cite{dougherty jackson kechris:1994}, two countably infinite non-smooth hyperfinite Borel equivalence relations are Borel
isomorphic if and only if the sets $EM_1(\Rk_1)$ and $EM_1(\Rk_2)$ of ergodic probability measures have the same cardinality
(see \cite{dougherty jackson kechris:1994} for definitions). Now the result follows from
Lemma \ref{invariant measures} and Theorem
\ref{Theorem_MeasuresErgodicMeasures}. The second statement is immediate from the first one as a special case. \hfill$\square$

\subsection{Aperiodic substitutions.} Let $\Ak$ denote a finite alphabet and $\Ak^+$ the set of all non-empty words over $\Ak$.
A map $\sigma:\, \Ak \to \Ak^+$ is called  a {\em substitution}. By concatenation, $\sigma$ is extended to the map $\sigma:\,\Ak^+
\to \Ak^+$. We define the {\em language} of the substitution $\sig$ as the set of all words which appear as factors of
$\sig^n(a),\ a\in \Ak,\ n\ge 1$. The {\em substitution dynamical system associated to $\sig$} is a pair $(X_\sig,T_\sig)$ where
$$
X_\sig = \{x\in \Ak^\Z:\, x[-n,n] \in \Lk(\sig)\ \mbox{for all}\ n\},
$$
and $T_\sig$ is the left shift on $\Ak^\Z$. The substitution $\sig$ is called {\em aperiodic} if the system $(X_\sig,T_\sig)$ has no
periodic points. Let $|w|$ denote the length of a word $w$. We will assume that
\begin{equation} \label{eq-long}
|\sig^n(a)| \to \infty,\ \mbox{as}\ n\to \infty,\ \mbox{for all}\ a\in \Ak.
\end{equation}
Such substitutions are sometimes called ``growing'', see \cite{Pansiot}. 
The substitution matrix is defined by $M_\sig = (m_{ab})$ where $m_{ab}$ is the number of letters $a$ occurring in $\sig(b)$. The
substitution $\sig$ is {\em primitive} if $M_\sig$ is primitive. It is well-known that primitive substitution dynamical systems are
minimal and uniquely ergodic, see \cite{queffelec:book}. 
There exists literature on non-primitive substitutions, including those for which (\ref{eq-long}) is violated, see e.g.\ 
\cite{AlSh,Pansiot,MS,Durand,DL},
mostly in the framework of combinatorics on words and theoretical computer science. However, the investigation of non-primitive, non-minimal,
{\em substitution dynamical systems}
has begun only recently \cite{yuasa,bezuglyi_kwiatkowski_medynets:substitutions}.

The connection between simple stationary Bratteli diagrams and primitive substitutions was first pointed out by Livshits \cite{Liv87,Liv88,VL},
and later clarified in \cite{forrest,durand_host_skau}. The extension to the aperiodic case was recently achieved in
\cite{bezuglyi_kwiatkowski_medynets:substitutions}).

Let $B = (V,E, \omega)$ be a stationary ordered Bratteli diagram. Choose a stationary labeling of $V_n$ by an alphabet $\mathcal A$: $V_n =
\{ v_n(a) : a \in \mathcal A\},\ n >0$. For $a\in \mathcal A$ we consider the vertex $v(a)$ in $V_n$, $n\ge 2$, and all the edges leading to it from $V_{n-1}$. 
These edges are coming from some vertices $v(a_1),\ldots,v(a_s)$, where we list them according to the $\om$-order. The map
$a \mapsto a_1\cdots a_s$ from $\mathcal A$ to
$\mathcal A^+$, does not depend on $n$ by stationarity  and determines a substitution
called the {\em substitution read on $B$}. The following result was proved
in  \cite{bezuglyi_kwiatkowski_medynets:substitutions}, extending \cite{forrest,durand_host_skau} from the primitive to the aperiodic case.

\begin{theorem}\label{substitutions_Bratteli_diagrams}  Let $B$ be a stationary $\omega$-ordered Bratteli diagram whose path space $X_B$ has no isolated points. 
Suppose $B$ admits an aperiodic Bratteli-Vershik system $(X_B,\varphi_\omega)$.  Then the system $(X_B,\varphi_\omega)$ is conjugate to an aperiodic substitution dynamical system (with substitution  read on $B$) if and only if no  restriction of
$\varphi_B$ to a minimal component is isomorphic to an odometer.

Conversely, assume that a substitution $\sigma$ satisfies (\ref{eq-long}).
Then the substitution
dynamical system is conjugate to the Vershik map of a stationary
ordered Bratteli diagram.
\end{theorem}
Actually, in the second part of the theorem, more general substitutions, those having a certain ``nesting property,'' are considered in
\cite{bezuglyi_kwiatkowski_medynets:substitutions}.
Let $\sigma : \mathcal A\to \mathcal A^+$ be an aperiodic
substitution satisfying (\ref{eq-long}). Denote by $B_\sigma$ the stationary Bratteli diagram
``read on the substitution''. This means that the substitution matrix $M_\sigma$ is  the
transpose to the incidence matrix of $B_\sigma$. It follows from
Theorem \ref{substitutions_Bratteli_diagrams} that there exists a
stationary ordered Bratteli diagram $B(X_\sigma, T_\sigma) = B$
whose Vershik map $\varphi_B$ is conjugate to $T_\sigma$. Thus, we
have two Bratteli diagrams associated to $(X_\sigma, T_\sigma)$.  It
follows from the results of Section \ref{stationary} that all
$T_\sigma$-invariant measures can be determined from the stationary
ordered Bratteli diagram $B$. We observe that the diagram $B$ may
have considerably more vertices than the diagram $B_\sigma$ (see
Example \ref{two_minimal_components} below). In fact, we can prove
the following statement.

\begin{theorem}\label{measures_and_substitutions}
There is a one-to-one correspondence $\Phi$ between the set of
ergodic $T_\sigma$-invariant  probability measures on the space
$X_\sigma$ and the set of ergodic $\mathcal R$-invariant probability
measures on the path space
$X_{B_\sigma}$ of the stationary diagram $B_\sigma$ defined by
substitution $\sigma$. The same statement holds for non-atomic infinite invariant measures.
\end{theorem}

{\it Proof}. Given an aperiodic substitution $\sigma$ defined on a
finite alphabet $\mathcal A$, construct the Bratteli diagram
$B_\sigma$ such that the substitution read on the diagram $B_\sig$
coincides with $\sigma$. Condition (\ref{eq-long}) implies that the tail equivalence relation $\Rk$ is
aperiodic. By rearranging the letters of
$\mathcal A$ we can assume that the incidence matrix of $B_\sigma$
has the form  (\ref{Frobenius Form}).
Obviously, $\sigma$ generates an ordering $\omega$ on $B_\sigma$.
Recall that the sets $X_{\max}(\omega)$ and $X_{\min}(\omega)$ are
finite. In
general, this ordering does not produce a Vershik map $\varphi $ on
the path space $X_{B_\sigma}$. However, it is clear that $\varphi$
is well defined at least for any infinite path from the
$\varphi$-invariant set $X_0 := X_{B_\sigma} \setminus
Orb_\varphi(X_{\max}(\omega)\cup X_{\min}(\omega))$.  Since $\Rk$ is
aperiodic, we see that $\varphi$ has no periodic points, and hence
every finite $\varphi$-invariant measure is non-atomic. This implies that
the dynamical systems $(X_{B_\sigma}, \mathcal R)$ and $(X_0,
\varphi)$ have the same set of ergodic invariant measures (both finite and infinite non-atomic).

Now consider the map $\pi : X_0 \to \mathcal A^{\Z}$ where $\pi(x) =
(\pi(x)_k)$ and $\pi(x)_k = a,\ a \in \mathcal A,$ if and only if
$\varphi^k(x)$ goes through the vertex $a\in V_1,\ k\in\Z$. Then
\begin{equation}\label{pi}
\pi\circ \varphi = T_\sigma\circ\pi.
\end{equation}

We will show that $\pi$ is injective on $X_0$. To do this, we use
the recognizability property proved in  \cite[Theorem
5.17]{bezuglyi_kwiatkowski_medynets:substitutions} for any aperiodic
substitution. It says that for any $\xi\in X_\sigma$ there exist a
unique $\eta\in X_\sigma$ and unique $i\in
\{0,1,...,|\sigma(\eta[0])| -1\}$ such that $\xi = T_\sig^i\sigma(\eta)$.

Take $\xi_1 \in X_\sigma$ and find $\xi_n \in X_\sigma$ and $i_n$
such that $\xi_n = T_\sig^{i_{n+1}}\sigma(\xi_{n+1}),\ n \in \N.$ In
other words, $\xi_n$ and $\xi_{n+1}$ are related as follows (this is
an illustrative example):

$$
\begin{array}{cc|c|c|c|c|c|c|c|c}\hline \ldots & \multicolumn{1}{|c|} {\xi_n[-3]} &
\xi_n[-2] & \xi_n[-1] & \xi_n[0] & \xi_n[1] & \xi_n[2] & \xi_n[3] & \xi_n[4]  & \ldots \\
\hline \ldots & \multicolumn{2}{|c|} {\sigma(\xi_{n+1}[-1])} &
\multicolumn{3}{|c|} {\sigma(\xi_{n+1}[0])} & \multicolumn{3}{|c|}
{\sigma(\xi_{n+1}[1])} & \ldots\\ \hline \\
\end{array}
$$

Thus, every $\xi_1 \in X_\sigma$ generates an infinite matrix whose
rows are defined by $\xi_n$ as in the diagram above. Denote by
$X'_\sigma$ the subset of $X_\sigma$ formed by those $\xi_1$ for
which the picture shown above is infinite to the left and to the right; in other words, the blocks $\xi_n[0]$ grow in both
directions.
It follows from \cite[Theorem
A.1]{bezuglyi_kwiatkowski_medynets:substitutions} that the
complement of $X'_\sigma$ in $X_\sigma$ is at most countable.

We will now define a map $\tau$ from $X'_\sigma$ to $X_{B_\sigma}$. Given $\xi_1 \in
X'_\sigma$, we will construct an infinite path $x\in X_{B_\sigma}$ by the following
rule: the path $x = (e_n)$ goes through the vertices $\xi_n[0]\in
V_n$ and  $e_n$ is the $i_n$-edge with respect to the order of
$r^{-1}(\xi_n[0])$ (recall that the diagram has single edges between
the top vertex $v_0$ and the vertices of the first level). It is not
hard to check that $\pi(X_0) = X'_\sigma$ and $\tau\circ\pi =
\mbox{id}$, proving that $\pi$ is injective. It follows from
(\ref{pi}) that $(X_0, \varphi)$ is topologically conjugate to
$(\pi(X_0), T_\sigma)$. Since $X_{B_\sigma}\setminus X_0$ and
$X_\sigma \setminus X_\sigma'$ are at most countable, and there are
no periodic points, the claim of the theorem follows.
 \hfill$\square$
\\

\myremark \label{rem-almost} 
If $\mu$ is an ergodic $\mathcal R$-invariant
measure and $\nu = \Phi(\mu)$ is an ergodic $T_\sigma$-invariant
measure, then the clopen values sets coincide for $\mu$ and $\nu$. 

Moreover, it follows from the proof of the last theorem that $(X_{B_\sig},\varphi,\mu)$ and $(X_\sig, T_\sig,\nu)$
are almost topologically, and even finitary, conjugate, see \cite{DeKea}.
\\

Now let $\sigma :\mathcal A \to \mathcal A^+$ be an aperiodic
substitution satisfying (\ref{eq-long}). Passing from $\sig$ to a power $\sig^k$ does not change the substitution dynamical system,
so we can assume without loss of generality that the substitution matrix satisfies (\ref{prop-I}).
\begin{corollary} \label{cor-subs}
Let $\sigma :\mathcal A \to \mathcal A^+$ be an aperiodic
substitution having the property (\ref{eq-long}), with a substitution matrix
$M_\sigma$ satisfying condition (\ref{prop-I}).
Then the set of ergodic
probability measures for $T_\sigma$ is in 1-to-1 correspondence with
the set of distinguished eigenvalues for  $M_\sigma$. The
substitution dynamical system is uniquely ergodic if and only if it
has a unique minimal component (i.e.\ $s=1$ in (\ref{Frobenius Form for A})), and its
Perron-Frobenius eigenvalue is the spectral radius of
$M_\sigma$. The set of infinite non-atomic ergodic invariant measures for $T_\sigma$, which are positive and finite on at
least one open set (depending on the measure), modulo a constant multiple, is in 1-to-1 correspondence with
the set of Perron-Frobenius eigenvectors of the diagonal blocks
of $M_\sigma$, which are not distinguished and not equal to the $1\times 1$ matrix $[1]$.
\end{corollary}

{\em Proof.} This is a combination of
Theorem~\ref{measures_and_substitutions},
Theorem~\ref{Theorem_MeasuresErgodicMeasures}, and
Theorem~\ref{th-infmeas}. The unique ergodicity claim follows from
the definition of distinguished eigenvalues. \hfill$\square$
\\

\myremark \label{rem-subst} 1. It seems plausible that aperiodicity
assumption in the above corollary may be dropped. Yuasa \cite{yuasa}
investigated almost minimal substitution dynamical systems, which
have a fixed point as the unique minimal component, for which he
obtained a similar statement. More precisely, he considered
$M_\sigma$ with two diagonal blocks: a $1\times 1$ block $[\ell]$,
with $\ell \ge 2$, which corresponds to the fixed point, and another
primitive block. Then the system is uniquely ergodic if and only if
$\ell$ is the spectral radius of $M_\sigma$; then there is also an
invariant $\sigma$-finite measure of full support. Earlier, the
special case of ``Cantor substitution'' $0 \to 000,\ 1\to 101$ was
considered by A. Fisher \cite{fisher}. Note that the results of
Yuasa are complementary to ours since we study the aperiodic case;
however, the general non-aperiodic case remains open.

2. F. Durand \cite{Durand} obtained ``a theorem of Cobham for
non-primitive substitutions'' for ``good'' substitutions. A substitution is ``good'' if there is a minimal component with the
PF eigenvalue equal to the spectral radius of $M_\sig$. If the substitution is aperiodic and the minimal component is unique,
then being ``good'' is equivalent to being uniquely ergodic.

\medskip

\myexample\label{two_minimal_components} Consider the substitution
$\sigma$ on the alphabet $\{a,b,c,d,1\}$:

$$
\sigma=\left\{\begin{array}{l} a \mapsto ab\\b \mapsto ba\\c \mapsto
cd\\d \mapsto dc\\1 \mapsto a111c\\ \end{array}\right.
$$
The substitution dynamical system $(X_\sigma, T_\sigma)$ has two
minimal components $C_1$ and $C_2$ and each of them is conjugate to
the Morse substitution system. The substitution matrix of $\sigma$ is
$$M(\sigma)=
\left(\begin{array}{rrrrr}
1 & 1 & 0 & 0 & 1 \\
1 & 1 & 0 & 0 & 0 \\
0 & 0 & 1 & 1 & 1 \\
0 & 0 & 1 & 1 & 0 \\
0 & 0 & 0 & 0 & 3
\end{array}\right)
$$

The Bratteli diagram read on the substitution is shown in Fig.\,3:

\unitlength = 0.5cm
 \begin{center}
 \begin{graph}(17,15)
\graphnodesize{0.4}
 \roundnode{V0}(8.5,14)
\roundnode{V11}(0.0,9) \roundnode{V12}(4.25,9)
\roundnode{V13}(8.5,9) \roundnode{V14}(12.75,9)
\roundnode{V15}(17.0,9)
\roundnode{V21}(0.0,5) \roundnode{V22}(4.25,5)
\roundnode{V23}(8.5,5) \roundnode{V24}(12.75,5)
\roundnode{V25}(17.0,5)
\roundnode{V31}(0.0,1) \roundnode{V32}(4.25,1)
\roundnode{V33}(8.5,1) \roundnode{V34}(12.75,1)
\roundnode{V35}(17.0,1)
\edge{V11}{V0} \edge{V12}{V0} \edge{V13}{V0} \edge{V14}{V0}
\edge{V15}{V0}
\edge{V21}{V11} \edge{V21}{V12} \edge{V22}{V11} \edge{V22}{V12}
\edge{V23}{V13} \edge{V23}{V14} \edge{V24}{V13} \edge{V24}{V14}
\edge{V25}{V11} \edge{V25}{V13} \bow{V25}{V15}{0.06}
\bow{V25}{V15}{-0.06} \edge{V25}{V15}
\edge{V31}{V21} \edge{V31}{V22} \edge{V32}{V21} \edge{V32}{V22}
\edge{V33}{V23} \edge{V33}{V24} \edge{V34}{V23} \edge{V34}{V24}
\edge{V35}{V21} \edge{V35}{V23} \bow{V35}{V25}{0.06}
\bow{V35}{V25}{-0.06} \edge{V35}{V25}
\end{graph}
\centerline{.\ .\ .\ .\ .\ .\ .\ .\ .\ .\ .\ .\ .\ .\ .\ .\ }
\vskip0.3cm Fig. 3
 \end{center}

For the matrix $A = M_\sigma$ the distinguished eigenvalues are 2,
2, and 3. The non-negative eigenvectors corresponding to these
eigenvalues are
$$
p_1 = (1/2, 1/2, 0, 0, 0)^T,\ p_2= (0, 0, 1/2, 1/2, 0)^T,\  p_3=
(2/9, 1/9, 2/9, 1/9, 1/3)^T
$$
By Theorems \ref{Theorem_MeasuresErgodicMeasures} and
\ref{measures_and_substitutions}, they define the ergodic
$T_\sigma$-invariant probability measures $\mu_1,\ \mu_2,$ and
$\mu_3$. The measures $\mu_1$ and $\mu_2$ are supported on the
minimal components $C_1$ and $C_2$. On the other hand, $\mu_3$ is
supported by $X \setminus(C_1\cup C_2)$.

The Bratteli diagram $B(X_\sigma, T_\sigma)$
constructed by the method used in
\cite{bezuglyi_kwiatkowski_medynets:substitutions} has considerably
more vertices than $B_\sigma$, see Fig\,4. (We note that the Bratteli-Vershik map on the diagram in Fig.\,3 is
not topologically conjugate to the substitution system, whereas the one shown in Fig.\,4 is.) 

\unitlength = 0.7cm
\begin{center}
\begin{graph}(18,15)
\graphnodesize{0.4}
 \roundnode{V0}(9.0,14)
\roundnode{V11}(0.0,9) \roundnode{V12}(0.818181818182,9)
\roundnode{V13}(1.63636363636,9) \roundnode{V14}(2.45454545455,9)
\roundnode{V15}(3.27272727273,9) \roundnode{V16}(4.09090909091,9)
\roundnode{V17}(4.90909090909,9) \roundnode{V18}(5.72727272727,9)
\roundnode{V19}(6.54545454545,9) \roundnode{V110}(7.36363636364,9)
\roundnode{V111}(8.18181818182,9) \roundnode{V112}(9.0,9)
\roundnode{V113}(9.81818181818,9) \roundnode{V114}(10.6363636364,9)
\roundnode{V115}(11.4545454545,9) \roundnode{V116}(12.2727272727,9)
\roundnode{V117}(13.0909090909,9) \roundnode{V118}(13.9090909091,9)
\roundnode{V119}(14.7272727273,9) \roundnode{V120}(15.5454545455,9)
\roundnode{V121}(16.3636363636,9) \roundnode{V122}(17.1818181818,9)
\roundnode{V123}(18.0,9)
\roundnode{V21}(0.0,5) \roundnode{V22}(0.818181818182,5)
\roundnode{V23}(1.63636363636,5) \roundnode{V24}(2.45454545455,5)
\roundnode{V25}(3.27272727273,5) \roundnode{V26}(4.09090909091,5)
\roundnode{V27}(4.90909090909,5) \roundnode{V28}(5.72727272727,5)
\roundnode{V29}(6.54545454545,5) \roundnode{V210}(7.36363636364,5)
\roundnode{V211}(8.18181818182,5) \roundnode{V212}(9.0,5)
\roundnode{V213}(9.81818181818,5) \roundnode{V214}(10.6363636364,5)
\roundnode{V215}(11.4545454545,5) \roundnode{V216}(12.2727272727,5)
\roundnode{V217}(13.0909090909,5) \roundnode{V218}(13.9090909091,5)
\roundnode{V219}(14.7272727273,5) \roundnode{V220}(15.5454545455,5)
\roundnode{V221}(16.3636363636,5) \roundnode{V222}(17.1818181818,5)
\roundnode{V223}(18.0,5)
\roundnode{V31}(0.0,1) \roundnode{V32}(0.818181818182,1)
\roundnode{V33}(1.63636363636,1) \roundnode{V34}(2.45454545455,1)
\roundnode{V35}(3.27272727273,1) \roundnode{V36}(4.09090909091,1)
\roundnode{V37}(4.90909090909,1) \roundnode{V38}(5.72727272727,1)
\roundnode{V39}(6.54545454545,1) \roundnode{V310}(7.36363636364,1)
\roundnode{V311}(8.18181818182,1) \roundnode{V312}(9.0,1)
\roundnode{V313}(9.81818181818,1) \roundnode{V314}(10.6363636364,1)
\roundnode{V315}(11.4545454545,1) \roundnode{V316}(12.2727272727,1)
\roundnode{V317}(13.0909090909,1) \roundnode{V318}(13.9090909091,1)
\roundnode{V319}(14.7272727273,1) \roundnode{V320}(15.5454545455,1)
\roundnode{V321}(16.3636363636,1) \roundnode{V322}(17.1818181818,1)
\roundnode{V323}(18.0,1)
\edge{V11}{V0} \edge{V12}{V0} \edge{V13}{V0} \edge{V14}{V0}
\edge{V15}{V0} \edge{V16}{V0} \edge{V17}{V0} \edge{V18}{V0}
\edge{V19}{V0} \edge{V110}{V0} \edge{V111}{V0} \edge{V112}{V0}
\edge{V113}{V0} \edge{V114}{V0} \edge{V115}{V0} \edge{V116}{V0}
\edge{V117}{V0} \edge{V118}{V0} \edge{V119}{V0} \edge{V120}{V0}
\edge{V121}{V0} \edge{V122}{V0} \edge{V123}{V0}
\edge{V21}{V12} \edge{V21}{V13} \edge{V22}{V14} \edge{V22}{V15}
\edge{V23}{V11} \edge{V23}{V16} \edge{V24}{V12} \edge{V24}{V15}
\edge{V25}{V14} \edge{V25}{V16} \edge{V26}{V11} \edge{V26}{V13}
\edge{V27}{V18} \edge{V27}{V19} \edge{V28}{V110} \edge{V28}{V111}
\edge{V29}{V17} \edge{V29}{V112} \edge{V210}{V18} \edge{V210}{V111}
\edge{V211}{V110} \edge{V211}{V112} \edge{V212}{V17}
\edge{V212}{V19} \edge{V213}{V14} \edge{V213}{V16} \edge{V214}{V113}
\edge{V214}{V114} \edge{V214}{V115} \edge{V214}{V116}
\edge{V214}{V117} \edge{V215}{V114} \edge{V215}{V115}
\edge{V215}{V116} \edge{V215}{V117} \edge{V215}{V118}
\edge{V216}{V114} \edge{V216}{V115} \edge{V216}{V116}
\edge{V216}{V118} \edge{V216}{V119} \edge{V217}{V112}
\edge{V217}{V120} \edge{V218}{V14} \edge{V218}{V121}
\edge{V219}{V110} \edge{V219}{V112} \edge{V220}{V18}
\edge{V220}{V122} \edge{V221}{V11} \edge{V221}{V123}
\edge{V222}{V112} \edge{V222}{V120} \edge{V223}{V11}
\edge{V223}{V121}
\edge{V31}{V22} \edge{V31}{V23} \edge{V32}{V24} \edge{V32}{V25}
\edge{V33}{V21} \edge{V33}{V26} \edge{V34}{V22} \edge{V34}{V25}
\edge{V35}{V24} \edge{V35}{V26} \edge{V36}{V21} \edge{V36}{V23}
\edge{V37}{V28} \edge{V37}{V29} \edge{V38}{V210} \edge{V38}{V211}
\edge{V39}{V27} \edge{V39}{V212} \edge{V310}{V28} \edge{V310}{V211}
\edge{V311}{V210} \edge{V311}{V212} \edge{V312}{V27}
\edge{V312}{V29} \edge{V313}{V24} \edge{V313}{V26} \edge{V314}{V213}
\edge{V314}{V214} \edge{V314}{V215} \edge{V314}{V216}
\edge{V314}{V217} \edge{V315}{V214} \edge{V315}{V215}
\edge{V315}{V216} \edge{V315}{V217} \edge{V315}{V218}
\edge{V316}{V214} \edge{V316}{V215} \edge{V316}{V216}
\edge{V316}{V218} \edge{V316}{V219} \edge{V317}{V212}
\edge{V317}{V220} \edge{V318}{V24} \edge{V318}{V221}
\edge{V319}{V210} \edge{V319}{V212} \edge{V320}{V28}
\edge{V320}{V222} \edge{V321}{V21} \edge{V321}{V223}
\edge{V322}{V212} \edge{V322}{V220} \edge{V323}{V21}
\edge{V323}{V221}
\end{graph}
\vskip0.3cm \centerline{.\ .\ .\ .\ .\ .\ .\ .\ .\ .\ .\ .\ .\ .\ .\
.\ .\ .\ .\ .\ .\ .\ .\ .\ .\ .\ .\  } \vskip0.7cm Fig. 4
\end{center}

\noindent \myexample Let $\sigma$ be the substitution defined on the
alphabet $A= \{a,b,1,2,3\}$ as follows:
$$
\sigma=\left\{\begin{array}{l} a \mapsto ab\\b \mapsto ba\\1 \mapsto
a111a\\2 \mapsto a22b\\3 \mapsto 133332\\ \end{array}\right.
$$
It is not hard to see that $(X_\sigma, T_\sigma)$ has a unique
minimal component $C_1$ defined by the subdiagram based on the symbols
$\{a,b\}$. The substitution matrix of $\sigma$ is

$$
M(\sigma)=\left(\begin{array}{rrrrr}
1 & 1 & 2 & 1 & 0 \\
1 & 1 & 0 & 1 & 0 \\
0 & 0 & 3 & 0 & 1 \\
0 & 0 & 0 & 2 & 1 \\
0 & 0 & 0 & 0 & 4
\end{array}\right).
$$
Positive eigenvalues of $M(\sigma)$ that have non-negative
eigenvectors are $2$ (found  from the $2 \times 2$  matrix
in the upper-left corner), $3$, and $4$. Notice that $m_{4,4} = 2$
is not a distinguished eigenvalue. The corresponding eigenvectors
are $p_1 = (1/2, 1/2, 0, 0, 0)^T$, $ p_2 =(1/2, 1/4, 3/8, 0, 0)^T$,
and $p_3= (1/4, 1/8, 1/4, 1/8, 1/4)^T$. In this case we have three
ergodic $T_\sigma$-invariant measures $\nu_1, \nu_2$, and $\nu_3$
built by the vectors $p_1, p_2$, and  $p_3$, respectively. The
measure $\nu_1$ is supported on the minimal component  $C_1$.

The Bratteli diagram read on the substitution $\sigma$ has the
following form:

\unitlength = 0.5cm
 \begin{center}
 \begin{graph}(17,15)
\graphnodesize{0.4}
 \roundnode{V0}(8.5,14)
\roundnode{V11}(0.0,9) \roundnode{V12}(4.25,9)
\roundnode{V13}(8.5,9) \roundnode{V14}(12.75,9)
\roundnode{V15}(17.0,9)
\roundnode{V21}(0.0,5) \roundnode{V22}(4.25,5)
\roundnode{V23}(8.5,5) \roundnode{V24}(12.75,5)
\roundnode{V25}(17.0,5)
\roundnode{V31}(0.0,1) \roundnode{V32}(4.25,1)
\roundnode{V33}(8.5,1) \roundnode{V34}(12.75,1)
\roundnode{V35}(17.0,1)
\edge{V11}{V0} \edge{V12}{V0} \edge{V13}{V0} \edge{V14}{V0}
\edge{V15}{V0}
\edge{V21}{V11} \edge{V21}{V12} \edge{V22}{V11} \edge{V22}{V12}
\bow{V23}{V11}{0.06} \bow{V23}{V11}{-0.06} \bow{V23}{V13}{0.06}
\bow{V23}{V13}{-0.06} \edge{V23}{V13} \edge{V24}{V11}
\edge{V24}{V12} \bow{V24}{V14}{0.06} \bow{V24}{V14}{-0.06}
\edge{V25}{V13} \edge{V25}{V14} \bow{V25}{V15}{0.06}
\bow{V25}{V15}{-0.06} \bow{V25}{V15}{0.03} \bow{V25}{V15}{-0.03}
\edge{V31}{V21} \edge{V31}{V22} \edge{V32}{V21} \edge{V32}{V22}
\bow{V33}{V21}{0.06} \bow{V33}{V21}{-0.06} \bow{V33}{V23}{0.06}
\bow{V33}{V23}{-0.06} \edge{V33}{V23} \edge{V34}{V21}
\edge{V34}{V22} \bow{V34}{V24}{0.06} \bow{V34}{V24}{-0.06}
\edge{V35}{V23} \edge{V35}{V24} \bow{V35}{V25}{0.06}
\bow{V35}{V25}{-0.06} \bow{V35}{V25}{0.03} \bow{V35}{V25}{-0.03}
\end{graph}
\centerline{.\ .\ .\ .\ .\ .\ .\ .\ .\ .\ .\ .\ .\ .\ .\ .\ .\ .\ .\
} \vskip0.3cm Fig. 5 \vskip0.2cm
 \end{center}

Using the methods of \cite{bezuglyi_kwiatkowski_medynets:substitutions} it can be computed that
the Bratteli diagram $B(X_\sigma, T_\sigma)$ has the
incidence matrix of size $26 \times 26$.

%
%

\section{Ergodic-theoretic properties}

In this section we study the dynamical systems on stationary
Bratteli diagrams $B$ from the ergodic-theoretic point of view. Our
results extend the work of A. Livshits \cite{Liv88} on minimal Vershik
maps and substitution systems, and our methods are rather similar to
those of Livshits. We should also note that for minimal substitution
systems absence of mixing was proved by Dekking and Keane \cite{DK},
and the characterization of eigenvalues was obtained by Host
\cite{Host}. For linearly recurrent systems eigenvalues were studied in \cite{CDHM,BDM}. Our results have some common features
with these papers, but they do not follow from them, since we are no longer in the minimal uniquely ergodic setting.

Let $B=(V,E)$ be a stationary Bratteli diagram. Throughout this section we assume that (\ref{prop-I}) holds, which can always be
achieved by telescoping and reordering the vertices.
As was noted before, we may consider the Vershik map $\varphi$ on
$X_B$ defined everywhere except the orbits of maximal and minimal
paths. This yields a measure preserving system
$(X_B,\varphi,\mu_\alpha)$ even when $\varphi$ cannot be extended to
a homeomorphism of $X_B$. Here $\mu_\alpha$ is the
ergodic invariant probability measure determined by a distinguished
vertex $\alpha$ from the reduced graph $R(A)$. With every such
$\alpha$ we associate the subdiagram $B_\alpha$ consisting of
vertices from the class $\mathcal E_\alpha$ with $A_\alpha \neq 0$
and edges connected them.

For $v\in V$,  $E(v_0,v)$ denotes the set of all finite paths from
$v_0$ to $v$. Clearly, the Vershik map is defined in the natural way
for any path $\ove\in E(v_0,v)$ if $\ove $ is not maximal. Then  for
every two such finite paths $\ove$ and $\ove'$ from $E(v_0, v)$
there exists an integer $Q=Q(\ove,\ove')$ such that $\varphi^Q(\ove)
= \ove'$. We denote by $[\ove]$ the cylinder subset of $X_B$
corresponding to a finite path $\ove$.

Since $B$ is stationary, we can (and will) identify the vertex set $V_n$, for
$n\ge 1$, with the set $\{1,\ldots,N\}$ to agree with the indexing
of rows and columns of the matrix $A$. We also consider a ``vertical shift''
map $S$ on the set of edges $E \setminus E_1$, so that $S(E_n) =
E_{n+1}$, i.e. if $s(e) = i \in V_{n-1}$ and $r(e) = j \in V_n$,
then $s(S(e)) = i \in V_{n},\ r(S(e)) = j \in V_{n+1}$. This
transformation naturally extends to finite paths starting from
vertices of level $n\ge 1$.

A pair of distinct finite paths $(\ovom,\ovom')$ with $s(\ovom) =
s(\ovom')$ and $r(\ovom) = r(\ovom')$ will be called a {\em diamond}
(there is a similar notion of ``graph diamond'' in symbolic
dynamics). The {\em length} of the diamond is the common length of
$\ovom,\ovom'$. Let $\Dk_\alpha$ denote the set of diamonds with
both $\ovom,\ovom'$ in $B_\alpha$.

Let $(\ov{\om},\ovom')$ be a diamond and let $\ovtau$ be any path
from $v_0$ to $s(\ovom)=s(\ovom')$. It is easy to see that
$$
P(\ovom,\ovom'):= Q(\ovtau\ovom,\ovtau\ovom')
$$
is independent of $\tau$. (Here and below $\ovtau\ovom$ denotes the
natural concatenation of finite paths.) Thus,
\begin{equation} \label{eq-ver}
\varphi^{P(\ovom,\ovom')}[\ovtau\ovom] = [\ovtau\ovom'],\ \
\mbox{for all}\ \tau\in E(v_0,s(\ovom)).
\end{equation}
Observe that if $(\ov{\om},\ovom')$ is a diamond, then
$(S^n(\ovom),S^n(\ovom'))$ is a diamond as well. Denote
$$
P_n(\ovom,\ovom'):= P(S^n(\ovom),S^n(\ovom')).
$$

\begin{lemma} \label{lem-nonmix}
Let $(\ovom,\ovom')$ be a diamond in $\Dk_\alpha$. Then there exists
$\delta>0$  such that for every finite path $\ov{e}$ with
$r(\ov{e})$ in a class which has access to $\alpha$ in $R(A)$
\begin{equation} \label{eq-nonmix}
\mu_\alpha(\varphi^{P_n(\ovom,\ovom')}[\ove] \cap [\ove]) \ge \delta
\mu_\alpha([\ove])
\end{equation}
for all $n$ sufficiently large.
\end{lemma}

Before proving the lemma we deduce the following corollary.

\begin{corollary} \label{cor-nonmix}
(i) The system $(X_B,\varphi,\mu_\alpha)$ is not strongly mixing.

(ii) For any aperiodic substitution $\sigma$ having the property (\ref{eq-long}), the substitution dynamical system $(X_\sig,T_\sig,\nu)$ is
not strongly mixing for any ergodic probability measure $\nu$. 
\end{corollary}

{\em Proof.} (i) We can find a diamond in $\Dk_\alpha$ and apply the
lemma. Let $P_n = P_n(\ovom,\ovom')$. If the system was mixing, we
would have for every finite path $\ove$
$$
\mu_\alpha(\varphi^{P_n}[\ove] \cap [\ove]) \to
\mu_\alpha([\ove])^2, \ \ \mbox{as}\ n\to \infty,
$$
since $|P_n|\to \infty$. Choosing $\ove$ long enough, we can make
sure that $\mu_\alpha([\ove])< \delta$ and get a contradiction with
(\ref{eq-nonmix}). 

(ii) This follows from part (i) and the results of subsection 5.2.
\hfill$\square$
\medskip

{\em Proof of Lemma~\ref{lem-nonmix}.} Without loss of generality,
we can assume that the diamond $(\ovom,\ovom')$ starts at level 1
(any other diamond is obtained by vertical shifting). We can also assume that
for every vertex $\alpha$ the matrix $F_\alpha$ is strictly
positive. Suppose $s(\ovom) = s(\ovom') = j \in V_1$ and $r(\ovom) =
r(\ovom') = j'\in V_k$, so the diamond has length $k-1$. Suppose
$|\ove|=m$ and $r(\ove) = i\in V_m$. For $n\ge m + 2$ denote by
$[\ove; S^n(\ovom)]$ the cylinder set consisting of paths from
$[\ove]$ which go along the path $S^n(\ovom)$ from levels $n+1$ to
$n+k$. It follows from (\ref{eq-ver}) that
$$
[\ove;  S^n(\ovom)] \subset [\ove]\cap \varphi^{-P_n}[\ove] =
\varphi^{-P_n}(\varphi^{P_n}[\ove] \cap [\ove]).
$$
Thus, the desired claim will follow if we prove that
$$
\mu_\alpha([\ove; S^n(\ovom)]) \ge \delta \mu_\alpha([\ove])
$$
where $\delta>0$ is independent of $n\ge m = |\ove|$. By Theorem
\ref{Theorem_MeasuresErgodicMeasures},
$$
\mu_\alpha([\ove])  = x_i \lam_\alpha^{-m+1},
$$
with $x_i>0$ since $i$ is in a class which has access to $\alpha$.
On the other hand,
$$
\mu_\alpha([\ove; S^n(\ovom)]) = x_{j'} \lam_\alpha^{-n-k} \cdot
N(i,j),
$$
where $N(i,j)$ is the number of paths from $i\in V_m$ to $j\in
V_{n+1}$. We have
$$
N(i,j) = (A^{n+1-m})_{i,j} \sim \lam_\alpha^{n+1-m}
$$
by (\ref{eq-asymp}), as $n\to \infty$. It follows that
$$
\frac{\mu_\alpha([\ove; S^n(\ovom])}{\mu_\alpha([\ove])} \sim
\frac{x_i}{x_{j'} \lam_\alpha^k},
$$
which is independent of $n$, as desired. \hfill$\square$

\begin{theorem} \label{th-eigen}
A complex number $\gam$ is an eigenvalue for the finite
measure-preserving system $(X_B,\varphi,\mu_\alpha)$ if and only if
for every diamond $(\ovom,\ovom')\in \Dk_\alpha$,
\begin{equation} \label{eq-eigen}
\gam^{P_n(\ovom,\ovom')} \to 1,\ \ \mbox{as}\ n\to \infty.
\end{equation}
Moreover, if the diagram satisfies condition (\ref{prop-I'}),
then for $\gam$ to be an
eigenvalue it is sufficient that (\ref{eq-eigen}) holds for all
diamonds of length $k \le 2$ in $\Dk_\alpha$.
\end{theorem}
{\em Proof of necessity.} There are several closely related
approaches; we follow \cite[Theorem 4.3]{SolTil}. Fix a diamond
$(\ovom,\ovom')\in \Dk_\alpha$. Let $f$ be a non-constant measurable
function on $X_B$ such that $f(\varphi x) = \gam f(x)$ for
$\mu_\alpha$-a.e.\ $x\in X_B$. By ergodicity, we can assume that
$|f|=1$ a.e. For any $\eps>0$ we can find a simple function $g =
\sum_{i\in \Ik} c_i \chi_{E_i}$ such that $E_i = [\ove_i]$ is the cylinder
set corresponding to a finite path $\ove_i$, $\{E_i: i \in \Ik\}$
forms a finite partition of $X_B$, and $\|f - g\|_1 < \eps$ where
$\| \cdot\|_1$ is the norm in  $L^1(X_B,\mu_\alpha)$. Suppose that
$n \ge \max\{|e_i| : i\in \Ik \} +2$, and let $P_n =
P_n(\ovom,\ovom')$. Consider the set
$$
A_n:= \bigcup_{i\in \Ik} (\varphi^{P_n}E _i \cap E_i).
$$
We claim that
$$
\mu_\alpha(A_n) = \sum_i \mu_\alpha(E_i \cap \varphi^{P_n} E_i) \ge
\sum_i \delta \mu_\alpha(E_i) = \delta
$$
where $\delta$ is the same as in Lemma \ref{lem-nonmix} Indeed, if
$\ove_i$ terminates in a vertex which has access to  $\alpha$ in
$R(A)$, then  (\ref{eq-nonmix}) applies to $E_i$, and otherwise,
$\mu_\alpha (E_i) =0$. We have
$$
\Jk  :=  \int_{A_n} |f(\varphi^{-P_n}x) - f(x) |\,d\mu_\alpha =
\mu_\alpha(A_n) |\gamma^{-P_n} - 1| \ge \delta |\gam^{P_n} - 1|,
$$
since $f$ is an eigenfunction. On the other hand,
\begin{eqnarray*}
\Jk &  \le & \int_{A_n} |f(\varphi^{-P_n}x) -
g(\varphi^{-P_n}x)|\,d\mu_\alpha + \int_{A_n} |g(\varphi^{-P_n}x) -
g(x)|\,d\mu_\alpha \\ & + & \int_{A_n} |g(x) - f(x)|\,d\mu_\alpha  <
2\eps.
\end{eqnarray*}
Indeed, the first and the third integrals are less than $\eps$ by
the choice of $g$, and the second integral is zero, since on
$\varphi^{P_n}E _i \cap E_i$ we have $g(x) = g(\varphi^{-P_n}x) =
c_i$. Combining the last two inequalities yields
$$
|\gam^{P_n} - 1| \le 2\eps/\delta,
$$
proving (\ref{eq-eigen}). \hfill$\square$
\medskip

{\em Proof of sufficiency in Theorem~\ref{th-eigen}.} By telescoping
the Bratteli diagram with respect to the levels $V_{nd+1}$ for some
$d\in \N$, we can assume that condition (\ref{prop-I'}) is satisfied. The new dynamical
system is measure-theoretically isomorphic to the original one, so
it has the same set of eigenvalues. Moreover, every diamond
$(\ovom,\ovom')$ of the telescoped diagram, with $s(\ovom) \in V_1$,
corresponds to a diamond $(\overline{\tau},\overline{\tau'})$ of the
original diagram, and
$P_n(\ovom,\ovom')=P_{nd}(\overline{\tau},\overline{\tau'})$. Thus,
if we prove sufficiency of (\ref{eq-eigen}) for the telescoped
diagram, the general case will follow as well.

We need two lemmas, which are rather standard. Their statements hold
for all diamonds, but we only need them for diamonds of length $k\le
2$.

\begin{lemma} \label{lem-recur}
Let $(\ovom,\ovom')$ be a diamond of length $k\le 2$, and let $P_n =
P_n(\ovom,\ovom')$. Then $P_n$ is a recurrent sequence satisfying
the recurrence relation of the characteristic polynomial of $A$.
More precisely, if $\det(zI-A) = z^N - d_1 z^{N-1} - \cdots - d_N$,
then
\begin{equation} \label{eq-recur}
P_{n+N} = d_1 P_{n+N-1} + \ldots + d_N P_n, \ \ \mbox{for all}\ n\in
\N.
\end{equation}
\end{lemma}

{\em Proof of the lemma.} 
First suppose that $(\ovom,\ovom')$ has
length 1. Without loss of generality, assume that $s(\ovom) =
s(\ovom') = j\in V_1$ and $r(\ovom) = r(\ovom') = j'\in V_2$. Then
$\ovom = (\om_1)$ and $\ovom' = (\om'_1)$, with $\om_1,\om'_1$ two
distinct edges between $j$ and $j'$. Let $\kappa$ and $\kappa'$ be
the positions of these edges in the ordered set $r^{-1}(j')$. Now it
is easy to see that
\begin{equation} \label{eq-diam1}
P_n(\ovom,\ovom') = (\kappa'-\kappa) h_{j}^{(n)},
\end{equation}
and (\ref{eq-recur}) holds, since it holds for all $h_w^{(n)} =
\sum_{i=1}^N (A^n)_{i,w}$. Here we use the Caley-Hamilton Theorem which says that matrices $A^n$, hence all
their matrix elements, satisfy (\ref{eq-recur}).

Now suppose that $(\ovom,\ovom')$ has length 2. Then, without loss
of generality, we can assume that $\ovom = (\om_1\om_2)$ and $\ovom'
= (\om'_1\om'_2)$ are distinct paths from $j\in V_1$ to $j'\in V_3$.
We can also assume that $i:=r(\om_1) \ne i':=r(\om'_1)$, otherwise,
the diamond decomposes into two diamonds of length 1. Suppose that
$\om_2 < \om'_2$ in the linear ordering $r^{-1}(j')$ (if not, switch
$\ovom$ and $\ovom'$; this results in changing the sign of $P_n$).
Now it is not hard to see that
\begin{equation} \label{eq-diam2}
P_n(\ovom,\ovom') =\sum_{e\in r^{-1}(i):\,\om_1 \le e}
h_{s(e)}^{(n)} + \sum_{e\in r^{-1}(j):\,\om_2 \le e < \om_2'}
h_{s(e)}^{(n+1)} + \sum_{e\in r^{-1}(i'):\, e < \om_1'}
h_{s(e)}^{(n)},
\end{equation}
which implies (\ref{eq-recur}), since, once again, it holds for each $h_w^{(n)}$.  \hfill$\square$

\begin{lemma} \label{lem-conver}
Let $(\ovom,\ovom')$ be a diamond of length $k\le 2$, and let $P_n =
P_n(\ovom,\ovom')$. If $\gam^{P_n}\to 1$, then the convergence is
geometric, that is, there exists $\rho \in (0,1)$ such that
$$
|\gam^{P_n} - 1| \le C\rho^n
$$
for some $C>0$.
\end{lemma}

{\em Proof of the lemma.} Let $\gam = e^{2\pi i \th}$, then
$\gam^{P_n}\to 1$ is equivalent to $P_n\th\to 0$ mod $\Z$. Let
$$
M = \left( \begin{array}{ccccc} 0 & 1 & 0 & \cdots & 0 \\
                                0 & 0 & 1 & \cdots & 0 \\
                                \cdots & \cdots & \cdots & \cdots & \cdots \\
                                 0 & 0 & 0 & \cdots & 1 \\
                                d_N & d_{N-1} & d_{N-2} & \cdots & d_1 \end{array} \right),
\ \ \ \ \bx_n = \left[ \begin{array}{c} P_n \th \\ P_{n+1}\th \\
\vdots \\ P_{n+N-1}\th  \end{array} \right].
$$
Then $\bx_{n+1} = M\bx_n$ and $P_n\th\to 0$ mod $\Z$  implies that
$M^n \bx_1 \to 0$ mod $\Z^N$, as $n\to \infty$. Now the claim
follows from \cite[Lemme 1]{Host}. \hfill$\square$

\medskip

{\em Continuation of the proof of sufficiency.} Choose an infinite
path $$x^{(0)}= (x^{(0)}_1, x^{(0)}_2,\ldots)$$ in $X_B$ such that
its every vertex lies in the class $\alpha$. In our notation, this
means $x^{(0)}$ is in $Y_\alpha$. It may be convenient to choose
$x^{(0)}$ to be ``constant'' (that is $S(x^{(0)}_n) =
x^{(0)}_{n+1}$), which is possible since in $B_\alpha$ all vertices
are connected by the property (I$'$), but this is not necessary. Let
$f(x^{(0)}) =1$. Now consider an arbitrary $x\in X_B$. If $x\not\in
X_\alpha$, then we can set $f(x)=1$ (or any other value), since the
set of such paths has zero $\mu_\alpha$ measure by
Lemma~\ref{lem-support}. Then we can suppose that the vertices of
$x$ lie in $B_\alpha$ for all levels $n\ge N$. Fix $n\ge N$ and
consider the vertex $r(x_n)\in V_n(B_\alpha)$. If $x^{(0)}$ passes
through $r(x_n)$, take $e_n:= x^{(0)}_{n+1}$, otherwise take $e_n$
to be any edge connecting $r(x_n)$ to $r(x^{(0)}_{n+1})$. Let
$$
f_n(x) = \gam^{-Q_n}, \ \mbox{where}\ Q_n = Q_n(x)\in \Z\  \mbox{is
such that}\ \varphi^{Q_n} (x[1,n]e_n) = x^{(0)}[1,n+1],
$$
which is well-defined (note that it may be negative). Finally, let
$$
f(x) = \lim_{n\to \infty} f_n(x).
$$
We are going to show that this limit exists. We claim that there
exist $C>0$ and $\rho\in (0,1)$ such that
\begin{equation} \label{eq-bound}
|f_{n+1}(x) - f_n(x)| = |\gam^{Q_{n+1} - Q_n}- 1|< C\rho^n,
\end{equation}
for $n$ sufficiently large. This will imply convergence of $f_n$ to
$f$. Observe that the pair of finite paths
$$
(\ovom^{(n)}, {\ovom'}^{(n)}) := (e_n x^{(0)}_{n+2}, x_{n+1}
e_{n+1})
$$
forms a diamond in $\Dk_\alpha$ of length 2, see Fig.\ 6.

\begin{center}
\begin{figure}[h]
$\begin{array}{c} \epsfxsize=1.5in
\epsffile{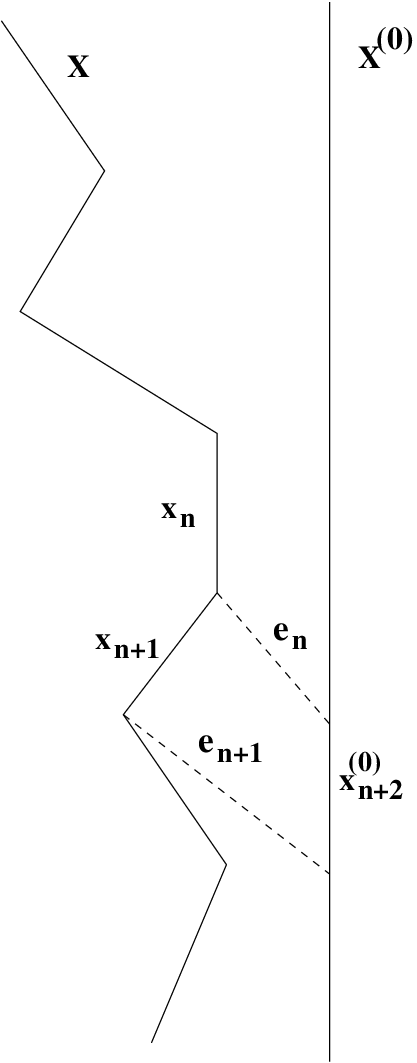} \\
[1.cm] \mbox{Fig.\ 6. The diamond $(\ovom^{(n)}, {\ovom'}^{(n)})$}
\end{array}$
\end{figure}
\end{center}

(It can be a ``degenerate diamond'' if the paths coincide, in which
case $Q_{n+1} = Q_n$ and there is nothing to prove.) In fact, there
is a diamond $(\ovom,\ovom')\in \Dk_\alpha$ starting at level 1 such
that $(\ovom^{(n)}, {\ovom'}^{(n)}) = (S^{n-1}(\ovom),
S^{n-1}(\ovom'))$. It remains to observe that $Q_{n+1} - Q_n =
P_{n-1}(\ovom,\ovom')$, so by Lemma~\ref{lem-conver} (keeping in
mind that there are finitely many possible diamonds of length 2) the
claim (\ref{eq-bound}) follows.

If we make a consistent choice of the edges $e_n$, it is clear that
this construction yields a measurable function $f$. In fact, $f_n$
are continuous on $X_\alpha$ and the convergence is uniform, so $f$
is continuous on $X_\alpha$ (however, we do not claim that $f$ has a
continuous extension as an eigenfunction to the entire $X_B$; this
need not be true).

It is easy to see that the definition of $f$ does not depend on the
choice of the edge $e_n$. This again follows from
Lemma~\ref{lem-conver}, since we get a diamond between levels $n$
and $n+1$ by choosing a different edge $e_n$. Finally, we claim that
$f$ is an eigenfunction. Since $x$ is a non-maximal path in
$X_\alpha$,  $\varphi(x)$ will only change the initial part of $x$
of certain length $k$. Take $n>k$ such that $x_n\in E(B_\alpha)$.
Then $(\varphi x)_n = x_n$ and we can choose the same edge $e_n$ in
the definition of $f(x)$ and $f(\varphi(x))$. It is clear that
$Q_n(\varphi(x)) = Q_n(x)-1$, hence $f_n(\varphi(x)) = \gam f_n(x)$,
and letting $n\to \infty$ we obtain $f(\varphi x) = \gam f(x)$, as
desired. \hfill$\square$
\medskip

\myremark \label{rem-cont} It is not hard to show by similar methods
that $\gam$ is an eigenvalue for the topological dynamical system
$(X_B,\varphi)$, with a continuous eigenfunction, if and only if
(\ref{eq-eigen}) holds for all diamonds in the diagram $B$.
Necessity is especially easy to see: if $(\ovom,\ovom')$ is a
diamond, then
$$
\varphi^{P_n(\ovom,\ovom')}[\ov{\tau}S^n(\ovom)]  =
[\ov{\tau}S^n(\ovom')],\ \ \mbox{for} \ \tau \in
E(v_0,s(S^m(\ovom)))
$$
by (\ref{eq-ver}). For any $x^{(n)}\in [\ov{\tau}S^n(\ovom)]$ we
obtain $\dist(x^{(n)}, \varphi^{P_n(\ovom,\ovom')}(x^{(n)}))\to 0$,
as $n\to \infty$, hence for a unimodular continuous eigenfunction
$f$ with eigenvalue $\gam$ we have by uniform continuity
$$
|f(x^{(n)}) - f(\varphi^{P_n(\ovom,\ovom')}(x^{(n)}))| = |1-
\gam^{P_n(\ovom,\ovom')}| \to 0, \ \ n\to \infty,
$$
as desired.

\medskip

Next we derive some consequences from Theorem~\ref{th-eigen}.

\begin{corollary} \label{cor-eigen}
Suppose that the Bratteli diagram satisfies condition (\ref{prop-I'}).
Then for $\gam = e^{2\pi
i\th}$ to be an eigenvalue of the measure-preserving system
$(X_B,\varphi,\mu_\alpha)$ it is sufficient that
\begin{equation} \label{eq-eigen2}
\th h_j^{(n)} \to 0\ {\rm mod}\ \Z,\ \ \mbox{as}\ n\to \infty,
\end{equation}
for every $j\in \Ek_\alpha$.
\end{corollary}

{\em Proof.} In view of Theorem~\ref{th-eigen}, this follows
from (\ref{eq-diam1}) and (\ref{eq-diam2}). 
\hfill$\square$

\medskip


Similarly to \cite{FMN}, 
it should be possible to determine the eigenvalues of
the system $(X_B,\varphi, \mu_\alpha)$, in an algebraic
way, and to obtain conditions for weak mixing. We do not pursue this
here, but restrict ourselves to a few illustrative examples. In
these examples, we specify the incidence matrix $F$ of the
stationary Bratteli diagram $B$. The matrix $F$ will be of size $2\times 2$ or $3\times 3$ and lower-triangular with at least one non-zero sub-diagonal entry
in each row
(except the first one, of course), so that $X_B$ will have a unique minimal component. Moreover, the  non-zero sub-diagonal entries will be all greater than one,
and we will define the linear order on $r^{-1}(v)$ in such a way that both the minimal and the maximal edge leading to $v$ come from another component
(except when $v$ is in the minimal component). Such an order produces a unique maximal and a unique minimal infinite path, which both lie in the minimal component.
So, the Bratteli-Vershik
homeomorphism $\varphi$  exists in each of the examples. We also assume that $h^{(1)} = (1,\ldots,1)^T$.
Recall that, in view of (\ref{heights}),
\begin{equation} \label{eq-hh1}
h^{(n+1)} = F^n h^{(1)}.
\end{equation}

\myexample \label{ex1} Let $F= \left( \begin{array}{cc} 2 & 0 \\ 2 & 3
\end{array} \right)$. 
There are two ergodic invariant probability measures on $X_B$: $\mu_1$, the
unique invariant measure on the minimal component, and $\mu_2$,
corresponding to the diagonal block $[3]$, which is fully supported. 

We will show that the system $(X_B,\varphi,\mu_2)$ has no non-trivial
eigenvalues, i.e.\ it is weakly mixing. 
An easy computation based on (\ref{eq-hh1}) yields
$$
h_1^{(n+1)} = 2^n,\ \ h_2^{(n+1)} = 3^{n+1} - 2^{n+1}.
$$
In the Bratteli diagram there exist two distinct edges $e_1,e_2$ leading
from the second vertex of $V_1$ to the second vertex of $V_2$, such that $e_2$ is
the immediate successor of $e_1$, producing a 
length-1 diamond $(\ovom,\ovom')\in \Dk_2$ 
with $\kappa'-\kappa=1$ in  (\ref{eq-diam1}). Thus, by (\ref{eq-diam1}), $P_n(\ovom,\ovom')= h_2^{(n)}= 3^{n} - 2^{n}$.
If
$\gam = e^{2\pi i \th}$ is an eigenvalue, then 
\begin{equation} \label{last1}
\th (3^{n} -
2^{n})\to 0\ \ \mbox{ mod}\  \Z,\ \ \mbox{as}\ n\to \infty,
\end{equation}
by Theorem~\ref{th-eigen}, and we claim that this implies $\gam=1$.
This can be shown by elementary considerations, but we refer the reader to a result of K\"orneyi \cite[Th.\,1]{Korn}, 
which we only partially quote here in a very special case.

\begin{theorem} \label{th-korn} (I. K\"orneyi) Let $\alpha_1,\ldots,\alpha_d$ be distinct integers, $|\alpha_j|\ge 1$ for $j\le d$, and $c_j\ne 0$ are
such that
$$
\sum_{j=1}^d c_j \alpha_j^n \to 0\ \ \mbox{\rm mod}\ \Z,\ \ \mbox{as}\ n\to \infty.
$$
Then $c_j \in \Q$ and $\sum_{j=1}^d c_j \alpha_j^n \in \Z$ for all $n$ sufficiently large.
\end{theorem}

In fact, in \cite{Korn} $\alpha_j$ are only assumed to be algebraic numbers, which is useful for determining
eigenvalues of Vershik maps in the general case. Returning to our example: by Theorem~\ref{th-korn} we infer from  (\ref{last1}) that
$\th(3^n-2^n)\in \Z$ for all $n$ sufficiently large, and it is elementary to check that then $\th$ is an integer, hence $\gam=1$.

The system $(X_B,\varphi,\mu_1)$ is isomorphic
to the 2-odometer, so it has pure discrete spectrum. As is well-known, and easily follows from Theorem~\ref{th-eigen},
$e^{2\pi i\th}$ is an eigenvalue for $(X_B,\varphi,\mu_1)$ if and only if $\th\cdot 2^n\to 0$\ mod $\Z$, as $n\to \infty$, that is,
$\th \in \Z[1/2]$. Notice,
however, that the eigenfunctions are not continuous on $X_B$ by
Remark~\ref{rem-cont}.

\medskip

The following examples show that the values of the
off-diagonal entries can affect the discrete spectrum. 

\myexample Let
$$
F = \left( \begin{array}{ccc} 5 & 0 & 0 \\ 2 & 3 & 0 \\ 0 &
2 & 25 \end{array} \right).
$$
We have a fully supported ergodic probability measure $\mu_3$ on $X_B$ corresponding to the eigenvalue $\lam_3 = 25$. Further,
$h^{(1)} = (1,1,1)^T = f_1 + (11/10)f_3$, where $f_1 = (1,1,-1/10)^T$ is the eigenvector of $F$ corresponding to
$\lam_1=5$, and $f_3 = (0,0,1)^T$ is the eigenvector corresponding to $\lam_3 = 25$. Then we obtain from (\ref{eq-hh1}):
$h^{(n+1)} = 5^n f_1 + (11/10)\cdot 25^n f_3$,
$$
h^{(n+1)} = (5^n,\ 5^n,\ (-5^n + 11\cdot 25^n)/10)^T.
$$
By Corollary~\ref{cor-eigen}, the set of eigenvalues for $(X_B,\phi_B,\mu_3)$
contains the set $\{\exp(2\pi p/5^n):\  n\ge 1,\ 1\le p \le 5^n-1\}$. 

\medskip

\myexample Let
$$
F = \left( \begin{array}{ccc} 5 & 0 & 0 \\ 4 & 3 & 0 \\ 0 &
2 & 25 \end{array} \right).
$$
The only difference from the previous example is the entry $F_{2,1}$; again we have 
a fully supported ergodic probability measure $\mu_3$ on $X_B$ corresponding to the eigenvalue $\lam_3 = 25$. However, in this case, the expression for
$h^{(1)}$ involves all three eigenvectors of $F$: $h^{(1)} = f_1 - f_2 + (61/55)f_3$, where
$$
f_1 = (1,2, -1/5)^T,\ \ \ f_2 = (0,1, -1/11)^T,\ \ \ f_3 = (0,0,1)^T.
$$
Thus,
$$ 
h^{(n+1)} = (5^n,\ 2\cdot 5^n - 3^n,\ -5^{n-1} +(1/11)3^n + (61/55)25^n)^T.
$$
We claim that the system $(X_B,\phi_B,\mu_3)$ is weakly mixing. The argument is similar to that of Example~\ref{ex1}. 
In the Bratteli diagram there exist two distinct edges $e_1,e_2$ leading from the third vertex of $V_1$ to the
third vertex of $V_2$, such that $e_2$ is
the immediate successor of $e_1$, producing a
length-1 diamond $(\ovom,\ovom')\in \Dk_3$
with $\kappa'-\kappa=1$ in  (\ref{eq-diam1}). Thus, by (\ref{eq-diam1}), $P_{n+1}(\ovom,\ovom')= h_3^{(n+1)}= K_n/55$ where
$K_n = -11\cdot 5^n + 5\cdot 3^n + 61\cdot 25^n$.
If
$\gam = e^{2\pi i \th}$ is an eigenvalue, then
$$
\th (K_n/55)
\to 0\ \ \mbox{ mod}\  \Z,\ \ \mbox{as}\ n\to \infty,
$$
by Theorem~\ref{th-eigen}.
By Theorem~\ref{th-korn}, we have 
\begin{equation} \label{eq-lll}
\th (K_n/55)\in \Z\ \ \mbox{ for all $n$ sufficiently large. }
\end{equation}
Let $ \th= p/q$, with $p,q$ mutually prime. Observe that $q$ is odd since $K_n$ is odd, and not divisible by $5$, because $K_n/55$ is not divisible by $5$. 
Next, note that $K_{n+1}-3K_n = 22\cdot 5^n \cdot (61\cdot 5^n - 1)$,
hence $(K_{n+1}-3K_n)/55 = 2 \cdot 5^{n-1} (61\cdot 5^n - 1)$. Any prime factor of $q$ must divide $61\cdot 5^n - 1$ (for all $n$ sufficiently large), hence
it is not $61$ and must divide $61(5^{n+1}-5^n)$ which does not contain any prime factors, other than 2, 5, and 61. We have proved that $\th$ is an integer,
hence $\gam=1$, as desired.
\\

{\it Acknowledgments}. The work was done during our mutual visits
to the University of Washington, University of Toru\'{n}, and
Institute for Low Temperature Physics. We are thankful to these
institutions for the hospitality and support.
We are grateful to the referee for comments and suggestions which helped improve the presentation, and to Alby Fisher for making his
preprints available prior to publication.

\end{document}